\documentclass[lettersize,journal]{IEEEtran}
\usepackage{amsmath,amssymb,amsfonts,amsthm}
\usepackage{cases}
\usepackage{algorithm}
\usepackage{algpseudocode}
\usepackage{array}
\usepackage[caption=false,font=normalsize,labelfont=sf,textfont=sf]{subfig}
\usepackage{textcomp}
\usepackage{stfloats}
\usepackage{url}
\usepackage{verbatim}
\usepackage{graphicx}
\usepackage{cite}
\usepackage{xcolor}
\usepackage{tikz}
\usetikzlibrary{shapes.geometric, arrows}
\usetikzlibrary{positioning}
\usepackage{booktabs}
\usepackage{multirow}
\usepackage{pgfplots} \pgfplotsset{compat=1.18}
\usepackage{siunitx,booktabs}
\usepackage{subfig}
\usepackage{hyperref}
\allowdisplaybreaks[3]

\newtheorem{theorem}{\bf Theorem}
\newtheorem{lemma}{\bf Lemma}
\newtheorem{definition}{\bf Definition}
\newtheorem{assumption}{\bf Assumption}
\newtheorem{problem}{\bf Problem}
\newtheorem{remark}{\bf Remark}

\begin{document}
\title{Reinforcement Learning for Discrete-time LQG Mean Field Social Control Problems with Unknown Dynamics}

\author{Hanfang Zhang, Bing-Chang Wang, \IEEEmembership{Senior Member, IEEE}, Shuo Chen
        % <-this % stops a space
\thanks{This work was supported in part by the National Natural Science Foundation of China under Grant No. 62192753, and Natural Science Foundation of Shandong Province for Distinguished Young Scholars (ZR2022JQ31).}% <-this % stops a space
%\thanks{Manuscript received April 19, 2021; revised August 16, 2021.}}
\thanks{H. Zhang and B.-C. Wang are with the School of Control Science and Engineering, Shandong University, Jinan, China (e-mail: \texttt{hanfangzhang@mail.sdu.edu.cn}, \texttt{bcwang@sdu.edu.cn}).}
\thanks{S. Chen is with the Beijing Institute for General Artificial Intelligence, Beijing, China (e-mail: \texttt{chenshuo@bigai.ai}).}
}
%\markboth{IEEE TRANSACTIONS ON SYSTEMS, MAN, AND CYBERNETICS: SYSTEMS,~Vol.~xx, No.~xx, August~2025}%
\markboth{Journal of \LaTeX\ Class Files,~Vol.~xx, No.~xx, August~2025}
{Zhang, Wang, and Chen: Reinforcement Learning for Discrete-time LQG Mean Field Social Control Problems}

%\IEEEpubid{0000--0000/00\$00.00~\copyright~2025 IEEE}

\maketitle

\noindent
This work has been submitted to the IEEE for possible publication. Copyright may be transferred without notice, after which this version may no longer be accessible.

\begin{abstract}
This paper studies the discrete-time linear-quadratic-Gaussian mean field (MF) social control problem in an infinite horizon, where the dynamics of all agents are unknown.
The objective is to design a reinforcement learning (RL) algorithm to approximate the decentralized asymptotic optimal social control in terms of two algebraic Riccati equations (AREs). 
In this problem, a coupling term is introduced into the system dynamics to capture the interactions among agents. This causes the equivalence between model-based and model-free methods to be invalid, which makes it difficult to directly apply traditional model-free algorithms.
First, under the assumptions of system stabilizability and detectability, a model-based policy iteration algorithm is proposed to approximate the stabilizing solution of the AREs. The algorithm is proven to be convergent in both cases of positive semidefinite and indefinite weight matrices.
Subsequently, by adopting the method of system transformation, a model-free RL algorithm is designed to solve for asymptotic optimal social control. 
During the iteration process, the updates are performed using
data collected from any two agents and the MF state.
Finally, a numerical case is provided to verify the effectiveness of the proposed algorithm.
\end{abstract}


\section{Introduction}
\IEEEPARstart{I}{n} recent years, mean field (MF) model has emerged as an important tool for modeling large-scale systems. The topic has been widely applied in various engineering fields, including unmanned aerial vehicles \cite{zhang2024large}, \cite{wang2022mean}, smart grids \cite{frihi2022hierarchical}, \cite{aziz2023computational}, intelligent urban rail transit \cite{li2023computation}, and epidemics \cite{bremaud2024mean}.
The MF game approach provides a critical theoretical framework for analyzing decentralized decision-making problems in large-scale multi-agent systems.
MF games originate from the parallel works of M. Huang et al. \cite{huang2006large, huang2007large} and of J. M. Lasry and P. L. Lions \cite{lasry2006jeux}, \cite{lasry2007mean}.
Inspired by these works, many fruitful results have been achieved (see, e.g., \cite{bensoussan2013mean, gomes2014mean, cardaliaguet2019master, wang2020mean}).

The core feature of MF games lies in the property that as the number of participants increases to a very large number, the influence of individual agents becomes negligible, while the impact of the population is significant.
Specifically, the interactions among individual agents are modeled through an MF term that represents the population aggregation effect,
thereby characterizing the high-dimensional game problem as a coupled system of forward-backward partial differential equations (the forward Kolmogorov-Fokker-Planck equation and the backward Hamilton-Jacobi-Bellman equation) \cite{lasry2007mean}.
As a classical type of MF models, linear-quadratic-Gaussian mean field (LQG-MF) has garnered particular attention due to its analytical tractability and practical approximation to physical systems \cite{huang2007large}, such as \cite{huang2010large, firoozi2022lqg, huang2019linear}.
The nonlinear MF games have the characteristic of their modeling generality (see, e.g. \cite{sen2016mean, fischer2017connection, lacker2020convergence, andrade2021stationary}). 

Social optima in MF models have attracted increasing attention. MF social control refers to that all agents cooperate to minimize a social cost as the sum of individual costs containing MF coupling term, which is generally regarded as a team decision-making problem.
For the early work, authors in \cite{huang2012social} investigated social optima in the LQG-MF control and provided an asymptotic team optimal solution, which was further extended to the case of mixed games in \cite{huang2016linear}. This model has also been applied to population growth modeling in \cite{yu2020social}. 
In the context of complex dynamic environments, \cite{wang2017social} adopted a parametric approach and state space augmentation to investigate the social optima of LQG-MF control models with Markov jump parameters.
\cite{wang2022robust} investigates the social optimality of MF control systems with unmodeled dynamics and applies it to analyzing opinion dynamics in social networks.
\cite{ma2023social} studied the MF social control problem with noisy output and designed a set of decentralized controllers by the variational method.
Furthermore, \cite{liang2024asymptotically} adopted the direct approach to investigate MF social control in a large-population system with heterogeneous agents.
For other aspects of MF control, readers may refer to \cite{mukaidani2025static}, \cite{barreiro2019semiexplicit} for nonlinear systems, \cite{huang2010optimal} for economic social welfare, \cite{salhab2018dynamic} for collective choice, and \cite{wang2024mean} for production output adjustment.

The aforementioned literature has made significant progress in the MF games and control problem. However, they all rely on the assumption that the system model is known.
In practical applications, complete system information is often difficult to obtain, and the system is susceptible to various external disturbances, which pose significant challenges to traditional control methods.
Reinforcement learning (RL) offers an effective approach to solving MF game and control problems with unknown system dynamics. 
Various methods have been developed, including fictitious play \cite{elie2020convergence}, Q-learning \cite{anahtarci2023learning}, and deep RL \cite{cui2021approximately}.
Additionally, \cite{uz2020reinforcement} proposed an actor-critic algorithm for RL in infinite-horizon non-stationary LQ-MF games.
\cite{lauriere2022scalable} gave two deep RL methods for dynamic MF games.
\cite{xu2023model} derived a set of decentralized strategies for continuous-time LQG-MF games based on the trajectory of a single agent and proposed a model-free method based on the Nash certainty equivalence-based strategy to solve $\epsilon$-Nash equilibria for a class of MF games.

Most of the work literature on RL algorithms focuses on non-cooperative MF games, while studies on cooperative MF social control remain relatively limited, which motivates us to conduct the present study.
For MF controls, \cite{gu2021mean} proposes an MF kernel-based Q-learning algorithm with a linear convergence rate.
In \cite{angiuli2022unified}, a unified two-timescale MF Q-learning algorithm was proposed, where the agent cannot observe the population's distribution. Both \cite{gu2021mean} and \cite{angiuli2022unified} model large-scale multi-agent systems as Markov decision process by defining the state/action space and transition probability function, and proposed Q-learning algorithms for MF control problems. 
In contrast, the works of \cite{wang2024online} and \cite{xu2025mean} describe the system dynamics through stochastic difference equations.
Specifically, \cite{wang2024online} developed an online value iteration algorithm for MF social control with ergodic cost functions.
Furthermore, \cite{xu2025mean} studied the continuous-time MF social optimization problem. They developed a novel model-free method that does not require any system matrices. Moreover, the algorithm improves computational efficiency by sampling the dataset from agents’ states and inputs.
Different from \cite{wang2024online} and \cite{xu2025mean}, we study MF social control in a discrete-time framework, which is more suitable for numerical algorithm design and computer implementation. On the one hand, we introduce coupling terms in the system dynamics to capture the interactions between agents. This makes it closer to reality, but difficult to apply traditional model-free algorithms directly. On the other hand, we relax the positive semidefinite requirement for the state and control weighting matrices as in \cite{xu2025mean}, and allow the weighting matrices to be not definite.

Motivated by the above literature on RL, this paper investigates a discrete-time infinite-horizon LQG-MF social control problem, in which the dynamics of agents are coupled by an MF coupling term and all system parameters are unknown.
We first propose a model-based policy iteration (PI) algorithm, and prove that the algorithm is convergent for different cases. 
When the weight matrices are positive semidefinite and the detectability condition is satisfied, the convergence of the iteration sequence can be ensured by applying the Lyapunov theorem.
When the weight matrices are indefinite, the Lyapunov theorem no longer holds. By analyzing the eigenvalues of the system matrix, we prove that the iteration sequence is monotonically decreasing and bounded below, which further implies the convergence of the algorithm.
In both cases, by selecting an appropriate initial value, the iterative sequence ultimately converges to the stabilizing solution of the coupled algebraic Riccati equations (AREs).
For the model-free control, a substantial challenge arises since the system parameters are fully unknown. Traditionally, deriving the solution of ARE benefits from the equivalence between model-based and model-free algorithms.
However, the MF coupling term in our case invalidates such equivalence and thus increases the complexity of solving the AREs. To overcome this challenge, we adopt a \emph{system transformation} approach that utilizes the state difference between two agents to eliminate the MF interactions, thereby restoring the equivalence between model-based and model-free approaches.
Through the system transformation, we design a model-free RL algorithm to solve for decentralized asymptotic optimal social control.
Notably, the algorithm uses a dataset sampled from the state trajectories and inputs associated with two agents and the MF coupling term.
By establishing the equivalence between model-free and model-based methods, we demonstrate the convergence of the RL algorithm.
Finally, the effectiveness of the algorithm is verified by a numerical example.

The contributions of this paper are summarized as follows:
\begin{itemize}
    \item For the discrete-time MF social control problems, a model-based PI algorithm is proposed, and its convergence is proven in different cases.
    In particular, when the state weighting matrix is indefinite, through eigenvalue analysis of relevant matrices, the algorithm is shown to converge to the unique stabilizing solution of the coupled AREs, which determines the feedback gain for MF asymptotic social control.
    \item For MF social control problems with unknown dynamics, a system transformation method is adopted to establish a data-driven iterative equation that eliminates the dependence of AREs on system matrices. Subsequently, we propose a model-free RL algorithm for obtaining the MF decentralized asymptotic optimal social control.
\end{itemize}

The remainder of this paper is organized as follows:
Section \uppercase\expandafter{\romannumeral 2} presents the MF social optimal control problem.
Section \uppercase\expandafter{\romannumeral 3} designs a model-based PI algorithm, which iteratively approximates the optimal solution of the AREs.
Section \uppercase\expandafter{\romannumeral 4} proposes a model-free RL algorithm to compute the optimal decentralized control set for MF social control with unknown system dynamics.
Section \uppercase\expandafter{\romannumeral 5} provides a numerical simulation to validate the effectiveness of the proposed algorithms.
Section \uppercase\expandafter{\romannumeral 6} concludes the paper and discusses future research directions.

%The following notation will be used throughout. 
A list of notations is presented as follows.
Let $\mathbb{R}^n$ be the $n$-dimensional Euclidean space and let $\mathbb{R}^{n \times m}$ be the set of real $n \times m$ matrices. Denote by $\mathbb{S}^n$ the set of symmetric matrices in $\mathbb{R}^{n \times n}$. For $A \in \mathbb{S}^n$ and $B \in \mathbb{R}^{n \times m}$, define $svec(A)\! := [a_{1\!1}, 2a_{1\!2}, \dots , 2a_{1\!n}, a_{2\!2}, 2a_{2\!3}, \dots, 2a_{n\!-\!1,n}, a_{n\!n}]^{T}\in \mathbb{R}^{n(n\!+\!1)/2}$ and $vec(B) := [b_{1}^{T}, b_{2}^{T}, \dots , b_{m}^{T}]^{T}\in \mathbb{R}^{nm}$, where $b_{1},\dots,b_{m}\in \mathbb{R}^{n}$. The identity matrix of size $n$ is written as $I_n$. The symbol $\otimes$ denotes the Kronecker product, and $\|\cdot\|$ denotes the Euclidean norm for vectors and the Frobenius norm for matrices. For a vector $z$ and a matrix $Q$, $\|z\|_{Q}^{2}=z^{T}Qz$. We use $Q>0$ ($Q\ge 0$) to mean that $Q$ is positive definite (positive semidefinite). The notation $\operatorname{tr}(Q)$ denotes the trace of $Q$.

\section{Problem Description}
Consider a large population system with $N$ agents, denoted as $\mathcal{A}=\{\mathcal{A}_{i}, 1 \leq i \leq N \}$, where $\mathcal{A}_{i}$ represents the $i$-th agent. The state of agent $i$ satisfies the following discrete-time linear stochastic difference equation
\begin{align}\label{xi_k+1}
	x_{i(k+1)}=Ax_{ik}+Gx_{k}^{(N)}+Bu_{ik}+Dw_{ik},    
\end{align}
where $x_{ik} \in \mathbb{R}^{n},u_{ik} \in \mathbb{R}^{m}$ are the state and control input for agent $i$, respectively. $x_{k}^{(N)}=\frac{1}{N}\sum_{i=1}^{N}x_{ik}$ is called the MF term. $\{w_{ik}, i=1,\dots, N\}$ is a sequence of independent random white noises with zero mean and variance $\sigma^{2}$. The coefficients $A, G, B$, and $D$ are assumed to be unknown deterministic matrices with appropriate dimensions.  
The cost function of the agent $i$ is given as
\begin{align}\label{Ji_u}
	J_{i}(u)=&\mathbb{E} \Big\{ \sum_{k=0}^{\infty}\gamma^{k}\big[\|x_{ik}-\Gamma x_{k}^{(N)}\|_{Q}^{2}+\|u_{ik}\|_{R}^{2}\big] \Big\},  
\end{align}
where $0<\gamma<1$ is the discount factor, $Q$, $R$ and $\Gamma$ are known, with $Q$ and $R$ being symmetric. The social cost for the system (\ref{xi_k+1}) and (\ref{Ji_u}) is defined as 
\begin{equation}
	J_{soc}= \sum_{i=1}^{N}J_{i}(u).
\end{equation}
The decentralized control set is defined below
\begin{align}
	\mathcal{U}_{d,i}=&\Big\{u_{i}|u_{ik}\in \mathbb{R}^{m} \ is \ \mathcal{F}_{i(k-1)} \ measurable,\cr
    &\sum_{k=0}^{\infty} \mathbb{E}\{u_{ik}^{T}u_{ik}\}< +\infty ,\ i=1,\dots,N\Big\},
\end{align}
where $\mathcal{F}_{ik}=\sigma \{x_{i0},w_{i0},w_{i1},\dots,w_{ik}\}$.

\begin{problem}\label{pro1}
Develop a data-driven method to find a set of decentralized control laws to optimize the social cost $J_{soc}$.    
\end{problem}
\begin{definition}\label{D1}
For any $i \in \{1, \dots,N\}$, system (\ref{xi_k+1}) is said to be $\gamma$-stabilizable 
\footnote{
A similar definition is $\rho$-stability. A system
$$ dx(t)=Ax(t)dt+Cx(t)dW(t)$$
is said to be $\rho$-stable (or $[A,C]$ is $\rho$-stable, for short) if $\mathbb{E}\int_{0}^{\infty} e^{-\rho t}\,\|x(t)\|^{2}\,dt < \infty$. See the work of Zhang et al \cite{zhang2008generalized} for further results about $\rho$-stability and $\rho$-stabilizability.
}
if there exists a feedback control law $u_{i}=K_{i}x_{i}$ such that for any initial state $x_{i0} \in \mathbb{R}^{n}$, the state $x_{ik}$ satisfies $\lim_{k \to \infty} \gamma^{k}\mathbb{E}\|x_{ik}\|^{2}=0$. 
In this case, $K_{i}$ is referred to as a $\gamma$-stabilizer of system (\ref{xi_k+1}).
\end{definition}

We make the following assumptions.
\begin{assumption}\label{A1}
$\{x_{i0}, i=1,\dots,N\}$ are mutually independent and also independent of $\{w_{i},i=1,\dots,N\}$. They have the same mathematical expectation and a finite second moment.
\end{assumption}

\begin{assumption}\label{A2}
	The system $(A, B)$ is $\gamma$-stabilizable, and the system $(A+G, B)$ is $\gamma$-stabilizable. 
\end{assumption}

\begin{assumption}\label{A3}
The system $(\sqrt{\gamma} A, \sqrt{Q})$ is detectable in the mean sense, and the system $(\sqrt{\gamma}(A+G), \sqrt{Q}(I-\Gamma))$ is detectable in the mean sense, $Q \geq 0$, $R > 0$. 
\end{assumption}

%\begin{remark}
%The system $(A,B)$ is $\gamma$-stabilizable, which is equivalent to the system $(A,B)$ is stabilizable in the mean sense. The same equivalence holds for the system $(A+G, B)$.   
%\end{remark}

\begin{remark}
The system $(A,B)$ is $\gamma$-stabilizable if and only if $(\sqrt{\gamma}\,A,B)$ is stabilizable in the mean square sense; the same holds for $(A+G,B)$.
\end{remark}

We define $K$ and $\bar{K}$ as follows, 
\begin{align}
    &K=\gamma(R+\gamma B^{T}PB)^{-1}B^{T}PA,\label{eq_K_gamma}\\
    &\bar{K}=\gamma(R+\gamma B^{T}\Pi B)^{-1}B^{T}\Pi (A+G)-K,\label{eq_bar_K_gamma}
\end{align}
where $P$ and $\Pi$ are solutions to the following AREs
\begin{numcases}{}
	P \!=\! \gamma A^{T}\!PA \!-\! \gamma^{2} A^{T}\!PB(R \!+\! \gamma B^{T}\!PB)^{-\!1}\!B^{T}\!PA \!+\! Q, \label{eq_P} \\ 
	\Pi \!=\! \begin{aligned}[t]
		& \gamma (A \!+\! G)^{T}\!\Pi(A \!+\! G) 
		\!-\! \gamma^{2} (A \!+\! G)^{T}\!\Pi B(R\!+\!\gamma B^{T}\!  \\
		& \times \Pi B)^{-1}\!B^{T}\!\Pi(A + G) + Q + Q_{\Gamma},
	\end{aligned} \label{eq_Pi}
\end{numcases}
with $Q_{\Gamma}\!=\!\Gamma^{T}\!Q\Gamma\!-\!Q\Gamma\!-\!\Gamma^{T}\!Q$. As an approximation to $x_{k}^{(N)}$, we obtain
\begin{align}\label{x_bar_AGBK}
    \bar{x}_{k+1}=\big(A+G-B(K+\bar{K})\big)\bar{x}_{k},
\end{align}
the decentralized control laws may be taken as
\begin{align}\label{eq_u_ik}
    \hat{u}_{ik}=-K\hat{x}_{ik}-\bar{K}\bar{x}_{k},
\end{align}
where
\begin{align}
    \hat{x}_{i(k+1)}=(A-BK)\hat{x}_{ik}+G\hat{x}_{k}^{(N)}-B\bar{K}\bar{x}_{k}+Dw_{ik}.
\end{align}
According to \cite{ma2023social}, the following result holds.
\begin{lemma}\label{lem01}
	Under Assumptions~\ref{A1}--\ref{A3}, for Problem \ref{pro1}, the set of decentralized control laws $\{\hat{u}_{1},\dots,\hat{u}_{N}\}$ given by (\ref{eq_u_ik}) has asymptotic social optimality, i.e.,
\begin{align}
	\big| \frac{1}{N}J_{soc}(\hat{u})-\frac{1}{N}\inf_{u \in \mathcal{U}_{c}}J_{soc}(u)\big|=O(\frac{1}{\sqrt{N}}),
\end{align}
where 
\begin{align}
    \nonumber
    \mathcal{U}_{c}\!=&\Big\{\!(u_{1},\dots,u_{N})|u_{ik}\!\in\! \mathbb{R}^{m} \ is\ \sigma\{\bigcup_{i=1}^{N} \! \mathcal{F}_{i(k\!-\!1)}\}  \ measurable,\\
    &\sum_{k=0}^{\infty} \mathbb{E}\{u_{ik}^{T}u_{ik}\}< +\infty \!\Big\}.
\end{align}
\end{lemma}
\begin{proof}
The proof follows the same arguments as those used in the proof of Theorem 4 in \cite{ma2023social}, and is therefore omitted for brevity.
\end{proof}

Let
\begin{align}
\bar{\mathcal{M}}=\left\{(P,\Pi)\!=\!(P^{T}\!,\Pi^{T})\mid \mathcal{H}(P)\! \geq\! 0, \bar{\mathcal{H}}(\Pi)\!\geq\! 0\right\},
\end{align}
where 
%\begin{small}
\begin{align*}
\mathcal{H}(P) &\!=\!
\begin{bmatrix}
	\gamma A^{T}PA\!-\!P\!+\!Q & \gamma A^{T}PB \\
	\gamma B^{T}PA         & R\!+\!\gamma B^{T}PB
\end{bmatrix},  \\
\bar{\mathcal{H}}(\Pi) &\!=\!
\begin{bmatrix}
	\gamma(A\!+\!G)^{T}\!\Pi(A\!+\!G)\!-\!\Pi\!+\!Q\!+\!Q_{\Gamma} \!&\! \gamma (A\!+\!G)^{T}\!\Pi B \\
	\gamma B^{T}\!\Pi(A\!+\!G)                              \!&\! R\!+\!\gamma B^{T}\!\Pi B
\end{bmatrix}\!.
\end{align*}
%\end{small}
\begin{assumption}\label{A4} \cite{li2003indefinite}
$\bar{\mathcal{M}} \neq \varnothing$ and has a nonempty interior $(\tilde{P},\tilde{\Pi})$ in the sense that $\mathcal{H}(\tilde{P})>0,\ \bar{\mathcal{H}}(\tilde{\Pi})>0$.
\end{assumption}
\begin{remark}
The discount factor $\gamma$ in the cost function (\ref{Ji_u}) guarantees that $J_{i}(u)$ is finite. assumption \ref{A3} can be weakened. Under assumption \ref{A4}, neither $R$ nor $Q$ is limited to be positive definite.
\end{remark}

For the case of indefinite weight $Q$, we have the following lemma on asymptotic optimality.
\begin{lemma}\label{lem02}
Let Assumptions \ref{A1}, \ref{A2} and \ref{A4} hold. Assume (\ref{eq_P})-(\ref{eq_Pi}) admit negative definite solutions $P^{-}$ and $\Pi^{-}$, respectively. Then, the decentralized control law (\ref{eq_u_ik}) achieves asymptotic social optimality. Furthermore, if the initial states $\{x_{i0}\}$ share the same variance, then the asymptotic average social optimum can be expressed as follows
\begin{align}\label{eq_J_soc_hat_u}
\lim_{N \to \infty}\frac{1}{N}J_{soc}(\hat{u})=\mathbb{E}[\|x_{i0}-\bar{x}_{0}\|_{P}^{2}+\|\bar{x}_{0}\|_{\Pi}^{2}]+q_{\infty},
\end{align}
where
\begin{align}
    q_{\infty}=\sum_{k=0}^{\infty}\gamma^{k}\,\sigma^{2}\,\mathrm{tr}\big(D^{T}(P+\Pi)D\big).
\end{align}
\end{lemma}
\begin{proof}
    See Appendix \ref{app:A}.
\end{proof}

Note that (\ref{eq_P}) and (\ref{eq_Pi}) are of vital importance to the design of the decentralized control.
However, it is challenging to solve the Riccati equations (\ref{eq_P}) and (\ref{eq_Pi}). We next develop a model-based MF social control design that leverages the system model to obtain stabilizing gains from the coupled AREs.

\section{Model-based MF social control design} \label{sec:Model-based-MF-social-control-design}
We approximate the decentralized control law (\ref{eq_u_ik}) by iteratively solving the Riccati equation, which is transformed into a Lyapunov equation. According to (\ref{eq_K_gamma}), equation (\ref{eq_P}) can be written as
\begin{align}\label{P_A}
	P=\gamma (A-BK)^{T}P(A-BK)+K^{T}RK+Q,  
\end{align}
For equation (\ref{P_A}), in order to approximate the sequence pairs $\{P,K\}$, a PI algorithm is presented as follows. 
We denote by $P_k$ the $k$-th iterative solution of the following Lyapunov equation
\begin{align}\label{eq_P_k} 
P_{k}= \gamma A_{k}^{T}P_{k}A_{k}+K_{k}^{T}RK_{k}+Q,
\end{align}
where $A_{k}=A-BK_{k}$, and update $K_{k+1}$ recursively by
\begin{align}\label{eq_K_k}
	K_{k+1}=\gamma (R+\gamma B^{T}P_{k}B)^{-1}B^{T}P_{k}A,  k=0,1,\cdots.
\end{align}
Based on Theorem 5 in \cite{dragan2020exact} and Proposition 3.2 in \cite{ni2015indefinite}, we present the following lemmas.
\begin{lemma}\label{lem03} \cite{dragan2020exact}
Suppose that $(A,\sqrt{Q})$ is detectable, there does not exist a non-zero symmetric matrix $\mathcal{Z}$, such that
\begin{align}
	\left\{
	\begin{aligned}
		&\mathcal{Z} - A^{T}\mathcal{Z}A=\lambda \mathcal{Z}, \ \lvert \lambda \rvert \geq 1,\\
		&\sqrt{Q}\mathcal{Z} = \mathbf{0}.   
	\end{aligned}
	\right.    
\end{align}
\end{lemma}
\vspace{0.5em}
\begin{lemma}\label{lem04}\cite{ni2015indefinite}
Under the assumption \ref{A2}, if $\bar{\mathcal{M}} \neq \varnothing$ and has a nonempty interior $(\tilde{P},\tilde{\Pi})$ in the sense that $\mathcal{H}(\tilde{P})>0,\ \bar{\mathcal{H}}(\tilde{\Pi})>0$, then the AREs (\ref{eq_P})-(\ref{eq_Pi}) admit a stabilizing solution.
\end{lemma}

The following theorem shows the convergence of the model-based PI method for two cases.

\begin{theorem}\label{thm01}
    Suppose assumption \ref{A1} holds, and $K_{0} \in \mathbb{R}^{m\times n}$ is a $\gamma$-stabilizer of system $(A,B)$. Let $P_{k}$ and $K_{k+1}$ be solutions to (\ref{eq_P_k})-(\ref{eq_K_k}), respectively. If either assumption \ref{A3} or assumption \ref{A4} is additionally satisfied, the following properties hold:
	\\ \noindent
	a)$\,$ For all $k \geq 0$, $\sqrt{\gamma}A_{k}$ is Schur;
	\\ \noindent
	b)$\,$ $P_{k} \geq P_{k+1} \geq P$;
	\\ \noindent
	c)$\,$ $\lim_{k \rightarrow \infty}P_{k}=P$ and $\lim_{k \rightarrow \infty}K_{k}=K$.
\end{theorem}
\begin{proof}
    See Appendix \ref{app:A}.
\end{proof}

For notational simplicity, from (\ref{eq_bar_K_gamma}) we obtain the following equation based on (\ref{eq_Pi}),
\begin{align}\label{Pi_A_G}
    \nonumber
	\Pi=&\gamma (A+G-BK-B\bar{K})^{T}\Pi(A+G-BK-B\bar{K})\\
    &+(K+\bar{K})^{T}R(K+\bar{K})+Q+Q_{\Gamma}.
\end{align}
Denote the $k$-th iteration equation $\Pi_{k}$ based on (\ref{Pi_A_G}), which corresponds to the policy evaluation equation.
\begin{align}\label{eq_Pi_k}
	\Pi_{k}=\gamma \bar{A}_{k}^{T}\!\Pi_{k}\bar{A}_{k}\!+\!(K\!+\!\bar{K}_{k})^{T}\!R(K\!+\!\bar{K}_{k})\!+\!Q\!+\!Q_{\Gamma},
\end{align}
where $\bar{A}_{k}\!\!=\!\!A\!+\!G\!-\!\!B\!K\!\!-\!\!B\!\bar{K}_{k}$ and update $\bar {K}_{k\!+\!1}$ recursively by
\begin{align}\label{eq_K_k_bar}
\!\!\!\!\bar K_{k\!+\!1}\!\!=\!\gamma (R\!+\!\gamma B^{T}\!\Pi_{k} B)^{\!-1}\!B^{T}\!\Pi_{k}(A\!+\!G)\!-\!K, \!k\!=\!0,\!1,\!\cdots\!.
\end{align}
Similar to Theorem \ref{thm01}, we can obtain the following convergence result.
\begin{theorem}\label{thm02}
    Suppose assumption \ref{A1} holds, and $\bar{K}_{0}+K \in \mathbb{R}^{m\times n}$ is a $\gamma$-stabilizer of system $(A+G,B)$. Let $\Pi_{k}$ and $\bar{K}_{k+1}$ be solutions to (\ref{eq_Pi_k})-(\ref{eq_K_k_bar}), respectively. If either assumption \ref{A3} or assumption \ref{A4} is additionally satisfied, the following properties hold:
	\\ \noindent
	a)$\,$ For all $k \geq 0$, $\sqrt{\gamma}\bar{A}_{k}$ is Schur;
	\\ \noindent
	b)$\,$ $\Pi_{k} \geq \Pi_{k+1} \geq \Pi$;
	\\ \noindent
	c)$\,$ $\lim_{k \rightarrow \infty}\Pi_{k}=\Pi$ and $\lim_{k \rightarrow \infty}\bar{K}_{k}=\bar{K}$.
\end{theorem}
\begin{proof}
    See Appendix \ref{app:B}.
\end{proof}

In the subsequent steps, based on the iterative equations (\ref{eq_P_k})-(\ref{eq_K_k}) and (\ref{eq_Pi_k})-(\ref{eq_K_k_bar}), along with the MF state dynamics (\ref{x_bar_AGBK}), we aim to further eliminate the dependence on the system matrix coefficients in the AREs.

\section{Model-free MF social control design}
%在这一小节中，我们设计了一种无模型算法近似分散式控制策略集\ref{eq_u_ik}.由于系统动态方程中耦合项的存在，基于模型方法与无模型方法的等价性被打破，我们采用系统变换的方法重新恢复了该等价性。接下来，通过强化学习技术，消除了AREs对系统参数的依赖，最后基于所得到的增益矩阵，我们计算了平均场状态近似。
In this subsection, we propose a model-free algorithm to approximate the decentralized control policy set (\ref{eq_u_ik}). Due to the MF coupling term in our case, which disrupts the equivalence between model-based and model-free methods, we restore this equivalence through a system transformation approach. Then, by employing RL techniques, we eliminate the dependence of the AREs on system parameters. Finally, using the obtained gain matrices, we compute an approximation of the MF state.

\subsection{Matrix approximation with unknown dynamics}
To proceed, we define error variables and average variables
\begin{align}
    &\Delta x_{k} =\mathbb{E}[x_{ik}-x_{jk}], \Delta u_{k}=\mathbb{E}[u_{ik}-u_{jk}],i \neq j,\\
    %&\textcolor{blue}{\delta_{k} =x_{ik}-x_{jk},\qquad \mu_{k}=u_{ik}-u_{jk},i \neq j, }\cr
    &\bar{x}_{k} =\mathbb{E}[\frac{1}{N}\sum_{i=1}^{N}x_{ik}], \bar{u}_{k} =\mathbb{E}[\frac{1}{N}\sum_{i=1}^{N}u_{ik}],
\end{align}
By equation (\ref{xi_k+1}), the system dynamics can be written as
\begin{align}
    \Delta x_{k+1}=&A\Delta x_{k}+B\Delta u_{k}\cr
    =&A_{k}\Delta x_{k}+B(\Delta u_{k}+K_{k}\Delta x_{k}),\label{delta_x_k+1} \\
    \bar{x}_{k+1}=&(A+G)\bar{x}_{k}+B\bar{u}_{k}\cr
    =&\bar{A}_{k}\bar{x}_{k}+B(\bar{u}_{k}+(K+\bar{K}_{k})\bar{x}_{k})\label{bar_x_k+1}.
\end{align}
Define the following quadratic function
\begin{align*}
    V_{1}(\Delta x_{k})=\Delta x_{k}^{T}P_{k}\Delta x_{k}.
\end{align*}
By equation (\ref{delta_x_k+1}), one has
\begin{align}\label{delta_K+1_k_1}
    &\mathbb{E}[\gamma \Delta x_{k+1}^{T}P_{k}\Delta x_{k+1}-\Delta x_{k}^{T}P_{k}\Delta x_{k}]\cr
    =&\mathbb{E}[\Delta x_{k}^{T}(\gamma A^{T}P_{k}A-P_{k})\Delta x_{k}+2\gamma \Delta u_{k}^{T}B^{T}P_{k}A\Delta x_{k}\cr
    &+\gamma \Delta u_{k}^{T}B^{T}P_{k}B\Delta u_{k}],
\end{align}    
which is equivalent to
\begin{align}\label{delta_K+1_k_2}
    &\mathbb{E}[\gamma\Delta x_{k+1}^{T}P_{k}\Delta x_{k+1}-\Delta x_{k}^{T}P_{k}\Delta x_{k}]\cr
    =&\mathbb{E}[\Delta x_{k}^{T}(\gamma A_{k}^{T}P_{k}A_{k}+2\gamma K_{k}^{T}B^{T}P_{k}A_{k}+\gamma K_{k}^{T}B^{T}P_{k}BK_{k}\cr
    &-P_{k})\Delta x_{k}+2\gamma \Delta u_{k}^{T}(B^{T}P_{k}A_{k}+B^{T}P_{k}BK_{k})\Delta x_{k}\cr
    &+\gamma \Delta u_{k}^{T}B^{T}P_{k}B\Delta u_{k}].
\end{align}      
By substituting equations (\ref{eq_K_k}) and (\ref{eq_P_k}) into equation (\ref{delta_K+1_k_2}) and then letting
\begin{align}
    &Q_{k}=K_{k}^{T}RK_{k}+Q,\\
    &\Lambda_{k}^{1}=\gamma B^{T}P_{k}B,\\
    &\mathcal{K}_{k+1}=(R+\Lambda_{k}^{1})K_{k+1},
\end{align}
we can get
\begin{align}\label{delta_K+1_k_3}
    &\mathbb{E}[\gamma\Delta x_{k+1}^{T}P_{k}\Delta x_{k+1}-\Delta x_{k}^{T}P_{k}\Delta x_{k}]\cr
    =&\mathbb{E}[-\Delta x_{k}^{T}Q_{k}\Delta x_{k}\!+\!2\Delta u_{k}^{T}\mathcal{K}_{k\!+\!1}\Delta x_{k}\!+\!2\Delta x_{k}^{T}K_{k}^{T}\mathcal{K}_{k\!+\!1}\Delta x_{k}\cr
    &+\Delta u_{k}^{T}\Lambda_{k}^{1}\Delta u_{k}-\Delta x_{k}^{T}K_{k}^{T}\Lambda_{k}^{1}K_{k}\Delta x_{k}].   
\end{align}     
In addition, by Kronecker product representation, we have
\begin{align*}
    &\Delta x_{k}^{T}Q_{k}\Delta x_{k}=(\Delta x_{k}^{T} \otimes \Delta x_{k}^{T})svec(Q_{k}),\cr
    &\Delta x_{k}^{T}K_{k}^{T}\mathcal{K}_{k+1}\Delta x_{k}=(\Delta x_{k}^{T} \otimes \Delta x_{k}^{T})(I_{n} \otimes K_{k}^{T}) vec(\mathcal{K}_{k+1}),\cr
    &\Delta x_{k}^{T}K_{k}^{T}\Lambda_{k}^{1}K_{k}\Delta x_{k}=(\Delta x_{k}^{T} \otimes \Delta x_{k}^{T})(K_{k}^{T} \otimes K_{k}^{T}) svec(\Lambda_{k}^{1}).
\end{align*}    
Let $l>0$ represent the number of sets of training data. The following definitions are made for the sake of convenience.
\begin{equation}\label{eq_matrix01}
\left\{
\begin{aligned}
&\mathcal{I}_{\Delta x\Delta x}\triangleq [\mathcal{I}_{\Delta x\Delta x}^{1},\mathcal{I}_{\Delta x\Delta x}^{2},\dots,\mathcal{I}_{\Delta x\Delta x}^{l-1}],\\
&\mathcal{I}_{\Delta x\Delta x}^{k}\triangleq \mathbb{E}[\Delta x_{k}^{T} \otimes \Delta x_{k}^{T}], \\
&\mathcal{I}_{\Delta x\Delta x}^{'}\triangleq [\mathcal{I}_{\Delta x\Delta x}^{2},\mathcal{I}_{\Delta x\Delta x}^{3},\dots,\mathcal{I}_{\Delta x\Delta x}^{l}],\\
&\mathcal{I}_{\Delta x\Delta u}\triangleq [\mathcal{I}_{\Delta x\Delta u}^{1},\mathcal{I}_{\Delta x\Delta u}^{2},\dots,\mathcal{I}_{\Delta x\Delta u}^{l-1}],\\
&\mathcal{I}_{\Delta x\Delta u}^{k}\triangleq \mathbb{E}[\Delta u_{k}^{T} \otimes \Delta x_{k}^{T}], \\
&\mathcal{I}_{\Delta u\Delta u}\triangleq [\mathcal{I}_{\Delta u\Delta u}^{1},\mathcal{I}_{\Delta u\Delta u}^{2},\dots,\mathcal{I}_{\Delta u\Delta u}^{l-1}],\\
&\mathcal{I}_{\Delta u\Delta u}^{k}\triangleq \mathbb{E}[\Delta u_{k}^{T} \otimes \Delta u_{k}^{T}], 
\end{aligned}
\right.
\end{equation}
where $\mathcal{I}_{\Delta x\Delta x}, \mathcal{I}_{\Delta x\Delta x}^{'} \!\in\! \mathbb{R}^{(l-1) \times \frac{n(n+1)}{2}}$, $\mathcal{I}_{\Delta x\Delta u} \!\in\! \mathbb{R}^{(l-1) \times nm}$,  $\mathcal{I}_{\Delta u\Delta u} \!\in\! \mathbb{R}^{(l-1) \times \frac{m(m+1)}{2}}$.
%进而，公式（delta_K+1_k_3）可以重新写为如下的线性方程的矩阵形式
Then we have
\begin{align}\label{eq_I_xx_Qk_1}
    \mathcal{I}_{\Delta x\Delta x}&svec(Q_{k})\!=\!(\mathcal{I}_{\Delta x\Delta x}\!-\!\gamma\mathcal{I}_{\Delta x\Delta x}^{'})svec(P_{k})\cr
    &+\!2\mathcal{I}_{\Delta x\Delta u}vec(\mathcal{K}_{k+1})\!+\!2\mathcal{I}_{\Delta x\Delta x}(I_{n}\!\otimes\! K_{k}^{T})vec(\mathcal{K}_{k+1})\cr
    &+\!\mathcal{I}_{\Delta u\Delta u}svec(\Lambda^{1}_{k})\!-\!\mathcal{I}_{\Delta x\Delta x}(K_{k}^{T}\!\otimes\! K_{k}^{T})svec(\Lambda^{1}_{k}).
\end{align}
Further, equation (\ref{eq_I_xx_Qk_1}) can be rewritten in the following matrix form of a linear equation,
\begin{align}
    \mathfrak{A}_{k}^{1}\left[
		\begin{array}{cccc}
			svec(P_{k})  \\
			vec(\mathcal{K}_{k+1})\\
			svec(\Lambda_{k}^{1})
		\end{array}\right]=\mathfrak{B}_{k}^{1},
\end{align}
where
\begin{align*}
   \mathfrak{A}_{k}^{1}=&\big[\mathcal{I}_{\Delta x\Delta x}-\gamma \mathcal{I}_{\Delta x\Delta x}^{'}, 2\mathcal{I}_{\Delta x\Delta u}+2\mathcal{I}_{\Delta x\Delta x}(I_{n}\otimes K_{k}^{T}),\\
   &\mathcal{I}_{\Delta u\Delta u}-\mathcal{I}_{\Delta x\Delta x}(K_{k}^{T}\otimes K_{k}^{T})\big],\cr
   \mathfrak{B}_{k}^{1}=&\mathcal{I}_{\Delta x\Delta x}svec(Q_{k}).
\end{align*}    
\begin{assumption}\label{A5}
$l-1 \geq \frac{n}{2}(n+1)+nm+\frac{m}{2}(m+1)$ and
\begin{align}\label{eq_rank_1}
\nonumber
rank(\mathfrak{I}^{1})&=
rank \begin{pmatrix}
\begin{bmatrix} \mathcal{I}_{\Delta x\Delta x}^{1},\mathcal{I}_{\Delta x\Delta x}^{2},\dots,\mathcal{I}_{\Delta x\Delta x}^{l-1} \\
\mathcal{I}_{\Delta x\Delta u}^{1},\mathcal{I}_{\Delta x\Delta u}^{2},\dots,\mathcal{I}_{\Delta x\Delta u}^{l-1}\\ 
\mathcal{I}_{\Delta u\Delta u}^{1},\mathcal{I}_{\Delta u\Delta u}^{2},\dots,\mathcal{I}_{\Delta u\Delta u}^{l-1}
\end{bmatrix}\end{pmatrix}\\
&=\frac{n}{2}(n+1)+nm+\frac{m}{2}(m+1).  
\end{align}
\end{assumption}
\vspace{0.8em}
\begin{theorem}\label{thm03}
 Suppose assumption \ref{A5} holds,  then the unknown sequence $\{P_{k},\mathcal{K}_{k+1},\Lambda_{k}^{1}\}_{1}^{\infty}$ can be solved using the following equation
    \begin{align}\label{eq_P_K_lam}
		\begin{bmatrix} svec(P_{k})  \\
			vec(\mathcal{K}_{k+1})\\
			svec(\Lambda_{k}^{1}) \end{bmatrix}
            =(\mathfrak{A}_{k}^{1T}\mathfrak{A}_{k}^{1})^{-1}\mathfrak{A}_{k}^{1T}\mathfrak{B}_{k}^{1},% k=1,2,\dots,
	\end{align}
and it satisfies the following expression\\
    a) \ $\lim_{k \rightarrow \infty}P_{k}=P$; \\
    b) \ $\lim_{k \rightarrow \infty}(R+\Lambda_{k}^{1})^{-1}\mathcal{K}_{k+1}=K$.
\end{theorem}
\begin{proof}
    See Appendix \ref{app:B}.
\end{proof}

Second, we continue to eliminate the system information in the iterative equations (\ref{eq_K_k_bar})-(\ref{eq_Pi_k}).
Define the following quadratic function
\begin{align*}
    V_{2}(\bar{x}_{k})=\bar{x}_{k}^{T}\Pi_{k}\bar{x}_{k}.
\end{align*}
By equation (\ref{bar_x_k+1}), one has
\begin{align}\label{bar_K+1_k_1}
    &\mathbb{E}[\gamma \bar{x}_{k+1}^{T}\Pi_{k}\bar{x}_{k+1}-\bar{x}_{k}^{T}\Pi_{k}\bar{x}_{k}]\cr
    =&\mathbb{E}[\bar{x}_{k}^{T}(\gamma (A+G)^{T}\Pi_{k}(A+G)-\Pi_{k})\bar{x}_{k}\cr
    &+2\gamma \bar{u}_{k}^{T}B^{T}\Pi_{k}(A+G)\bar{x}_{k}+\gamma \bar{u}_{k}^{T}B^{T}\Pi_{k}B\bar{u}_{k}],
\end{align}    
which is equivalent to
\begin{align}\label{bar_K+1_k_2}
    &\mathbb{E}[\gamma \bar{x}_{k+1}^{T}\Pi_{k}\bar{x}_{k+1}-\bar{x}_{k}^{T}\Pi_{k}\bar{x}_{k}]\cr
    \!\!=&\mathbb{E}[\bar{x}_{k}^{T}\!\big(\gamma \bar{A}_{k}^{T}\!\Pi_{k}\bar{A}_{k}\!+\!2\gamma (K\!+\!\bar{K}_{k})^{T}\!B^{T}\!\Pi_{k}\bar{A}_{k}\!+\!\gamma (K\!+\!\bar{K}_{k})^{T}\cr
    &\times B^{T}\Pi_{k}B(K+\bar{K}_{k})-\Pi_{k}\big)\bar{x}_{k}+2\gamma \bar{u}_{k}^{T}\big(B^{T}\Pi_{k}\bar{A}_{k}\cr
    &+B^{T}\Pi_{k}B(K+\bar{K}_{k})\big)\bar{x}_{k}+\gamma \bar{u}_{k}^{T}B^{T}\Pi_{k}B\bar{u}_{k}].
\end{align}      
By substituting equations (\ref{eq_K_k_bar}) and (\ref{eq_Pi_k}) into equation (\ref{bar_K+1_k_2}) and then letting
\begin{align}
    &\bar{Q}_{k}=(K+\bar{K}_{k})^{T}R(K+\bar{K}_{k})+Q+Q_{\Gamma},\\
    &\Lambda_{k}^{2}=\gamma B^{T}\Pi_{k}B,\\
    &\mathcal{\bar{K}}_{k+1}=(R+\Lambda_{k}^{2})(K+\bar{K}_{k+1}),
\end{align}
we can get
\begin{align}\label{bar_K+1_k_3}
    &\mathbb{E}[\gamma \bar{x}_{k+1}^{T}\Pi_{k}\bar{x}_{k+1}-\bar{x}_{k}^{T}\Pi_{k}\bar{x}_{k}]\cr
    =&\mathbb{E}[-\bar{x}_{k}^{T}\bar{Q}_{k}\bar{x}_{k}+2\bar{u}_{k}^{T}\mathcal{\bar{K}}_{k+1}\bar{x}_{k}+2\bar{x}_{k}^{T}(K+\bar{K}_{k})^{T}\mathcal{\bar{K}}_{k+1}\bar{x}_{k}\cr
    &+\bar{u}_{k}^{T}\Lambda_{k}^{2}\bar{u}_{k}-\bar{x}_{k}^{T}(K+\bar{K}_{k})^{T}\Lambda_{k}^{2}(K+\bar{K}_{k})\bar{x}_{k}].
\end{align}     
The following definition is given by the properties of the Kronecker product,
\begin{equation}\label{eq_matrix02}
\left\{
\begin{aligned}
&\mathcal{I}_{\bar{x}\bar{x}}\triangleq [\mathcal{I}_{\bar{x}\bar{x}}^{1},\mathcal{I}_{\bar{x}\bar{x}}^{2},\dots,\mathcal{I}_{\bar{x}\bar{x}}^{l-1}],\mathcal{I}_{\bar{x}\bar{x}}^{k}\triangleq \mathbb{E}[\bar{x}_{k}^{T} \otimes \bar{x}_{k}^{T}], \cr
&\mathcal{I}_{\bar{x}\bar{x}}^{'}\triangleq [\mathcal{I}_{\bar{x}\bar{x}}^{2},\mathcal{I}_{\bar{x}\bar{x}}^{3},\dots,\mathcal{I}_{\bar{x}\bar{x}}^{l}],\cr
&\mathcal{I}_{\bar{x}\bar{u}}\triangleq [\mathcal{I}_{\bar{x}\bar{u}}^{1},\mathcal{I}_{\bar{x}\bar{u}}^{2},\dots,\mathcal{I}_{\bar{x}\bar{u}}^{l-1}],\mathcal{I}_{\bar{x}\bar{u}}^{k}\triangleq \mathbb{E}[\bar{u}_{k}^{T} \otimes \bar{x}_{k}^{T}], \cr
&\mathcal{I}_{\bar{u}\bar{u}}\triangleq [\mathcal{I}_{\bar{u}\bar{u}}^{1},\mathcal{I}_{\bar{u}\bar{u}}^{2},\dots,\mathcal{I}_{\bar{u}\bar{u}}^{l-1}],\mathcal{I}_{\bar{u}\bar{u}}^{k}\triangleq \mathbb{E}[\bar{u}_{k}^{T} \otimes \bar{u}_{k}^{T}], 
\end{aligned}
\right.
\end{equation}
where $\mathcal{I}_{\bar{x}\bar{x}},\ \mathcal{I}_{\bar{x}\bar{x}}^{'} \in \mathbb{R}^{(l-1)\times \frac{n(n+1)}{2}},\  \mathcal{I}_{\bar{x}\bar{u}} \in \mathbb{R}^{(l-1)\times nm}, \ \mathcal{I}_{\bar{u}\bar{u}} \in \mathbb{R}^{(l-1)\times \frac{m(m+1)}{2}}$.
Then we have
\begin{align}\label{eq_I_xx_Qk_1}
    &\mathcal{I}_{\bar{x}\bar{x}}svec(\bar{Q}_{k})\!=\!(\mathcal{I}_{\bar{x}\bar{x}}\!-\!\gamma\mathcal{I}_{\bar{x}\bar{x}}^{'})svec(\Pi_{k})\cr
    &+\!2\mathcal{I}_{\bar{x}\bar{u}}vec(\bar{\mathcal{K}}_{k+1})\!+\!2\mathcal{I}_{\bar{x}\bar{x}}(I_{n}\!\otimes\! (K\!+\!\bar{K}_{k})^{T})vec(\bar{\mathcal{K}}_{k+1})\ \ \ \cr
    &+\!\mathcal{I}_{\bar{u}\bar{u}}svec(\Lambda^{2}_{k})\!-\!\!\mathcal{I}_{\bar{x}\bar{x}}((K\!\!+\!\bar{K}_{k})^{T}\!\!\otimes\! (K\!\!+\!\bar{K}_{k})\!^{T})svec(\Lambda^{2}_{k}).\ \ \ 
\end{align}
Further, equation (\ref{bar_K+1_k_3}) can be rewritten in the following matrix form of a linear equation,
\begin{align}
    \mathfrak{A}_{k}^{2}\left[
		\begin{array}{cccc}
			svec(\Pi          _{k})  \\
			vec(\mathcal{\bar{K}}_{k+1})\\
			svec(\Lambda_{k}^{2})
		\end{array}\right]=\mathfrak{B}_{k}^{2},
\end{align}
where
\begin{align*}
   \mathfrak{A}_{k}^{2}=&[\mathcal{I}_{\bar{x}\bar{x}}-\gamma \mathcal{I}_{\bar{x}\bar{x}}^{'}, 2\mathcal{I}_{\bar{x}\bar{u}}+2\mathcal{I}_{\bar{x}\bar{x}}(I_{n}\otimes (K\!+\!
   \bar{K}_{k})^{T}), \mathcal{I}_{\bar{u}\bar{u}}\\
   &-\mathcal{I}_{\bar{x}\bar{x}}((K\!+\!\bar{K}_{k})^{T}\!\otimes\! (K\!+\!\bar{K}_{k})^{T})] ,\cr
   \mathfrak{B}_{k}^{2}=&\mathcal{I}_{\bar{x}\bar{x}}svec(\bar{Q}_{k}).
\end{align*}    
\begin{assumption}\label{A6}
$l-1 \geq \frac{n}{2}(n+1)+nm+\frac{m}{2}(m+1)$ and
\begin{align}\label{eq_rank_2}
rank(\mathfrak{I}^{2})&=
 rank \begin{pmatrix}
\begin{bmatrix} \mathcal{I}_{\bar{x}\bar{x}}^{1},\mathcal{I}_{\bar{x}\bar{x}}^{2},\dots,\mathcal{I}_{\bar{x}\bar{x}}^{l-1} \\
\mathcal{I}_{\bar{x}\bar{u}}^{1},\mathcal{I}_{\bar{x}\bar{u}}^{2},\dots,\mathcal{I}_{\bar{x}\bar{u}}^{l-1}\\ 
\mathcal{I}_{\bar{u}\bar{u}}^{1},\mathcal{I}_{\bar{u}\bar{u}}^{2},\dots,\mathcal{I}_{\bar{u}\bar{u}}^{l-1}
\end{bmatrix}\end{pmatrix}\cr
&=\frac{n}{2}(n+1)+nm+\frac{m}{2}(m+1). 
\end{align}
\end{assumption}
\begin{theorem}\label{thm04}
%如果假设成立，那么未知的序列对{P，K，Lambda}可以由下列公式求解得到
 Suppose assumption \ref{A6} holds,  then the unknown sequence $\{\Pi_{k},\bar{\mathcal{K}}_{k+1},\Lambda_{k}^{2}\}_{1}^{\infty}$ can be solved using the following equation
    \begin{align}\label{eq_Pi_K_lam}
		\begin{bmatrix} svec(\Pi_{k})  \\
			vec(\mathcal{\bar{K}}_{k+1})\\
			svec(\Lambda_{k}^{2}) \end{bmatrix}
            =(\mathfrak{A}_{k}^{2T}\mathfrak{A}_{k}^{2})^{-1}\mathfrak{A}_{k}^{2T}\mathfrak{B}_{k}^{2},% k=1,2,\dots,
	\end{align}
and it satisfies the following expression\\
    a) \ $\lim_{k \rightarrow \infty}\Pi_{k}=\Pi$; \\
    b) \ $\lim_{k \rightarrow \infty}\big((R+\Lambda_{k}^{2})^{-1}\mathcal{\bar{K}}_{k+1}-K\big)=\bar{K}$.
\end{theorem}
\begin{proof}
To prove that equation (\ref{eq_Pi_K_lam}) has a unique solution, we need to show that the matrix $\mathfrak{A}_{k}^{2}$ is of full column rank. The convergence result follows Theorem \ref{thm02}. We next show that $\mathfrak{A}_{k}^{2}$ is of full column rank.

Assume that there exists a vector $\mathcal{T}=[svec(\mathcal{T}_{1}),vec(\mathcal{T}_{2}),svec(\mathcal{T}_{3})]^{T}$, such that
\begin{align}
    \mathfrak{A}_{k}^{2}\mathcal{T}=\boldsymbol{0},
\end{align}
where $\mathcal{T}_{1} \!\in\! \mathbb{R}^{\frac{n(n+1)}{2}},\mathcal{T}_{2} \!\in\! \mathbb{R}^{m\times n},\mathcal{T}_{3} \!\in\! \mathbb{R}^{\frac{m(m+1)}{2}}$. Then we have
\begin{align}\label{eq_prf_3}
    \nonumber
    (\mathcal{I}_{\bar{x}\bar{x}}&-\gamma \mathcal{I}_{\bar{x}\bar{x}}^{'})svec(\mathcal{T}_{1})+\mathcal{I}_{\bar{x}\bar{x}}svec\big((K+\bar{K}_{k})^{T}\mathcal{T}_{2}\\
    \nonumber
    &+\mathcal{T}_{2}^{T}(K+\bar{K}_{k})-(K+\bar{K}_{k})^{T}\mathcal{T}_{3}(K+\bar{K}_{k})\big)\\  
    &+2\mathcal{I}_{\bar{x}\bar{u}}vec(\mathcal{T}_{2})+\mathcal{I}_{\bar{u}\bar{u}}svec(\mathcal{T}_{3})=\boldsymbol{0},
\end{align}   
According to the equation (\ref{bar_K+1_k_1}), it gives
\begin{align}\label{eq_prf_4}
    \nonumber
    &(\gamma \mathcal{I}_{\bar{x}\bar{x}}^{'}-\mathcal{I}_{\bar{x}\bar{x}})svec(\mathcal{T}_{1})\\
    \nonumber
    =&\mathcal{I}_{\bar{x}\bar{x}}svec\big(\gamma (A+G)^{T}\mathcal{T}_{1}(A+G)-\mathcal{T}_{1}\big)\\
    &+2\mathcal{I}_{\bar{x}\bar{u}}vec\big(B^{T}\mathcal{T}_{1}(A+G)\big)+\mathcal{I}_{\bar{u}\bar{u}}svec(B^{T}\mathcal{T}_{1}B).
\end{align} 
Combining equation (\ref{eq_prf_3}) and (\ref{eq_prf_4}), we can get
\begin{align}
  \mathcal{I}_{\bar{x}\bar{x}} svec(\mathcal{V}_{1})+2\mathcal{I}_{\bar{x}\bar{u}} vec(\mathcal{V}_{2})+\mathcal{I}_{\bar{u}\bar{u}} svec(\mathcal{V}_{3})=\boldsymbol{0},  
\end{align}
where
\begin{align*}
    \mathcal{V}_{1}=&\gamma (A+G)^{T}\mathcal{T}_{1}(A+G)-\mathcal{T}_{1}-(K+\bar{K}_{k})^{T}\mathcal{T}_{2}\cr
    &
    -\mathcal{T}_{2}^{T}(K+\bar{K}_{k})+(K+\bar{K}_{k})^{T}\mathcal{T}_{3}(K+\bar{K}_{k}),\cr
    \mathcal{V}_{2}=&\gamma B^{T}\mathcal{T}_{1}(A+G)-\mathcal{T}_{2},\cr
    \mathcal{V}_{3}=&\gamma B^{T}\mathcal{T}_{1}B-\mathcal{T}_{3}.
\end{align*}
Based on the rank condition (\ref{eq_rank_2}), we can derive that $\mathcal{V}_{1}=\mathcal{V}_{2}=\mathcal{V}_{3}=\boldsymbol{0}$. Then we have
\begin{align}
    \mathcal{T}_{1}-\gamma \bar{A}_{k}^{T}\mathcal{T}_{1}\bar{A}_{k}=\boldsymbol{0},
\end{align}
%由于前面引理2证明Ak是舒尔的，所以可以得到S1=0，进而得到S2=S3=0.
Since Theorem \ref{thm02} previously proved that $\sqrt{\gamma} \bar{A}_{k}$ is Schur. According to \cite{hewer1971iterative}, we can conclude that $\mathcal{T}_{1}=0$, which in turn implies that $\mathcal{T}_{2}=\mathcal{T}_{3}=0$.
%由此，S=0,这与假设S不等于0矛盾，所以Ak是满秩的。
Thus, we have $\mathcal{T}=\boldsymbol{0}$,and so $\mathfrak{A}_{k}^{2}$ has full column rank.
\end{proof}

\begin{figure*}[htpb]
\centering
\resizebox{\textwidth}{!}{
\begin{tikzpicture}
    \node [rectangle, draw=none, fill=blue!5, rounded corners, minimum width=24.4cm, minimum height=2.3cm] (data_box) at (7.5,-2) {};
    \node [text=blue!60, font=\bfseries, anchor=west, xshift=-0.6cm, yshift=0cm] at (data_box.west){\parbox{3.5cm}{\centering \large Data \\Collection}};
    \node [rectangle, draw=none, fill=pink!20, rounded corners, minimum width=9.4cm, minimum height=6.5cm] (pk_box) at (0,-6.85) {};
    \node [text=pink!80!red, font=\bfseries, anchor=west, xshift=-0.2cm, yshift=-1.8cm] at (pk_box.west){\parbox{3.5cm}{\centering \large Approximation \\of $(P, K, \Lambda^{1})$}};
    \node [rectangle, draw=none, fill=purple!10, rounded corners, minimum width=9.4cm, minimum height=6.5cm] (pik_box) at (15 ,-6.85) {};
    \node [text=purple!60, font=\bfseries, anchor=east, xshift=0.2cm, yshift=-1.7cm] at (pik_box.east){\parbox{3.5cm}{\centering \large Approximation \\of $(\Pi, \bar{K}, \Lambda^{2})$}};

    \node (start) [rectangle, rounded corners, draw=black, fill=red!10, font=\large]{Start};
    \node (init) [rectangle, right=2cm of start, draw=black, fill=none, font=\large] {\parbox{9.5cm}{Choose a $\gamma$-stabilizer $K_{0}$ of system $(A,B)$, set $k=0$, discount factor $\gamma$ and convergence criterion $\epsilon$.}};
    \node (data01) [rectangle, below right=0.8cm and -2cm of init, draw=black, fill=none, font=\large] {\parbox{8cm}{Employ $u_{il}=-K_{0}x_{il}+\xi_{il}$, where $i=1,2$, to collect datasets $\mathcal{D}_{1},\ \mathcal{D}_{2}$. $l$ represents the iteration index, and $\xi_{il}$ is the exploration noise.}};
    \node (data02) [rectangle, left= 0.6cm of data01, draw=black, fill=none, font=\large] {\parbox{4.8cm}{Calculate matrices $\mathcal{I}_{\Delta x\Delta x}$, $\mathcal{I}_{\Delta x\Delta x}^{'}$, $\mathcal{I}_{\Delta x\Delta u}$, $\mathcal{I}_{\Delta u\Delta u}$, $\mathcal{I}_{\bar{x}\bar{x}}$, $\mathcal{I}_{\bar{x}\bar{x}}^{'}$, $\mathcal{I}_{\bar{x}\bar{u}}$, $\mathcal{I}_{\bar{u}\bar{u}}$.}};
    \node (data03) [diamond, left=0.6cm of data02, draw=black, fill=none, aspect=4.4, font=\large] {\parbox{4cm} {\centering
    $\mathfrak{I}^1$ and $\mathfrak{I}^2$ full rank?}};
    \node (PK01) [rectangle, below=0.8cm of data03, draw=black, fill=none, font=\large] {\parbox{7cm}{\centering Solve $\{P_{k},\mathcal{K}_{k+1},\Lambda_{k}^{1}\}$ by
    $
		\begin{bmatrix} svec(P_{k})  \\
			vec(\mathcal{K}_{k+1})\\
			svec(\Lambda_{k}^{1}) 
        \end{bmatrix}
        =(\mathfrak{A}_{k}^{1T}\mathfrak{A}_{k}^{1})^{-1}\mathfrak{A}_{k}^{1T}\mathfrak{B}_{k}^{1},
	$}};
    \node (PK02) [diamond, below=0.6cm of PK01, draw=black, fill=none, aspect=4, font=\large]{\parbox{4cm} {\centering
    $\| K_{k+1}-K_{k} \| > \epsilon$ ?}};
    \node (PK00) [rectangle, above left=0.2cm and 1.5cm of PK02, draw=black, fill=none, font=\large]{\parbox{1.6cm}{\centering $k=k+1$
    }};
    \node (PK03) [rectangle, below=0.8cm of PK02, draw=black, fill=none, font=\large]{\parbox{7cm}{\centering $\hat{P}=P_{k}$,\ $\hat{K}=(R+\Lambda_{k}^{1})^{-1}\mathcal{K}_{k+1}$
    }};
    \node (PiK) [rectangle, above right=1.5cm and 0.45cm of PK03, draw=black, fill=none, font=\large] {\parbox{4.8cm}{\centering Choose a $\gamma$-stabilizer $\bar{K}_{0}+\hat{K}$ of system $(A+G,B)$ and set $k=0$}};
    \node (PiK01) [rectangle, above right=0.1cm and 0.65cm of PiK, draw=black, fill=none, font=\large] {\parbox{7cm}{\centering Solve $\{\Pi_{k},\mathcal{\bar{K}}_{k+1},\Lambda_{k}^{2}\}$ by
    $
		\begin{bmatrix} vec(\Pi_{k})  \\
			vec(\mathcal{\bar{K}}_{k+1})\\
			vec(\Lambda_{k}^{2}) 
        \end{bmatrix}
        =(\mathfrak{A}_{k}^{2T}\mathfrak{A}_{k}^{2})^{-1}\mathfrak{A}_{k}^{2T}\mathfrak{B}_{k}^{2},% k=1,2,\dots,
	$}};
    \node (PiK02) [diamond, below=0.6cm of PiK01, draw=black, fill=none, aspect=4, font=\large]{\parbox{4cm} {\centering
    $\| \bar{K}_{k+1}-\bar{K}_{k} \| > \epsilon$ ?}};
    \node (PiK00) [rectangle, above right=0.2cm and 1.6cm of PiK02, draw=black, fill=none, font=\large]{\parbox{1.6cm}{\centering $k=k+1$
    }};
    \node (PiK03) [rectangle, below=0.8cm of PiK02, draw=black, fill=none, font=\large]{\parbox{7cm}{\centering $\hat{\Pi}=\Pi_{k}$,\ $\hat{\bar{K}}=(R+\Lambda_{k}^{2})^{-1}\mathcal{\bar{K}}_{k+1}-\hat{K}$
    }};
    \node (result) [rectangle, below left=0.5cm and -1.8cm of PiK03, draw=black, fill=none, font=\large]{\parbox{9.5cm}{\centering Apply $u_{ik}=-\hat{K}x_{ik}-\hat{\bar{K}}x_{k}^{(N)}$, where $i=1,2,...,N$.
    }};
    \node (end) [rectangle, rounded corners, left=2cm of result, draw=black, fill=red!10, font=\large] {End};   
    % Draw arrows
    \draw [thick,->] (start.east) -- (init.west);
    \draw [thick,->] (init.east) -- ++(2.12,0)  -- ++(0,-0.5) --  (data01.north);
    \draw [thick,->] (data01.west) -- (data02.east);
    \draw [thick,->] (data02.west) -- (data03.east);
    \draw [thick,->] (data03.south) -- node[midway, right] {Yes} (PK01.north);
    \draw [thick,->] (data03.north) --++ (0,0.3)  -- ++(6.5,0) node[pos=0.2, above] {No} |-++ (0,-0.2) -- (data02.north);
    \draw [thick,->] (PK01) --(PK02);
    \draw [thick,->] (PK02) -- node[midway, right] {Yes} (PK03);
    \draw [thick,->] (PK02.west) --++ (-0.7,0) node[pos=0.3, above] {No}  --++ (0,0.3) -- (PK00.south); 
    \draw [thick,->] (PK00.north) --++ (0,1.35) -- (PK01.west); 
    \draw [thick,->] (PK03.east) -- ++(3,0)-- ++(0,1)-- (PiK.south);
    \draw [thick,->] (PiK.north) -- ++(0,1.2)-- ++(1,0)-- (PiK01.west);
    \draw [thick,->] (PiK01) -- (PiK02);
    \draw [thick,->] (PiK02) -- node[midway, right] {Yes} (PiK03);
    \draw [thick,->] (PiK02.east) --++ (0.8,0) node[pos=0.3, above] {No} --++ (0,0.5)  -- (PiK00.south); % 连接到 k = k+1 框的底部
    \draw [thick,->] (PiK00.north) --++ (0,1.37) -- (PiK01.east); 
    \draw [thick,->] (PiK03.south) -- ++(0,-0.9)-- ++(-1.5,0) -- (result.east);
    \draw [thick,->] (result.west) -- (end.east);    
\end{tikzpicture}
}
\caption{Algorithm logic diagram}
\label{fig_logic}
\end{figure*}

\subsection{Data-driven MF social control algorithm design}
Herein, we are in the position to present the data-driven MF social control algorithm.
This algorithm eliminates the dependency on the system matrix in the model-based AREs discussed in Section 3. 
In the case of unknown system dynamics, the feedback gain matrices can be solved by updating them through the iterative optimization equation based on data-driven, and then the optimal control strategies can be obtained. 
The sampling dataset is collected from the states and inputs of any two agents, along with the relevant data of the MF state.
\begin{algorithm}[H]
    \caption{Data-driven model-free MF social control}
    \label{alg_01}
    \begin{algorithmic}[1]
    \Statex \hspace{-1.5em}\textbf{Input1:} \hspace{0.5em} Choose a $\gamma$-stabilizer $K_{0}$ of system $(A,B)$ and set convergence criterion $\epsilon$.
    \Statex \hspace{-1.5em}\textbf{Data:} \hspace{1.5em} Execute $u_{il}=-K_{0}x_{il}+\xi_{il}, \ i=1,2$, where $\xi_{ik}$ is the exploration noise and collect data $\mathcal{D}_{1}$, $\mathcal{D}_{2}$. 
    \Statex \hspace{3.5em}\textbf{repeat}
    \Statex \hspace{3.5em}\qquad Calculate matrices (\ref{eq_matrix01}), (\ref{eq_matrix02}).
    \Statex \hspace{3.5em}\textbf{until} Rank conditions (\ref{eq_rank_1}) and (\ref{eq_rank_2}) are satisfied.
    \Statex \hspace{-1.5em}\textbf{Output1:} $\hat{P}=P_{k}$,\ $\hat{K}=(R+\Lambda_{k}^{1})^{-1}\mathcal{K}_{k+1}$
    \Statex \hspace{3.5em}\textbf{while} $\| K_{k+1}-K_{k} \| > \epsilon$ \textbf{do}
    \Statex \hspace{3.5em}\qquad Solve $\{P_{k},\mathcal{K}_{k+1},\Lambda_{k}^{1}\}$ by (\ref{eq_P_K_lam}),
    \Statex \hspace{3.5em}\qquad $k=k+1$
    \Statex \hspace{3.5em}\textbf{end while}
    \Statex \hspace{-1.5em}\textbf{Input2:} \hspace{0.5em} Choose a $\gamma$-stabilizer $\bar{K}_{0}+\hat{K}$ of system $(A+G,B)$.
    \Statex \hspace{-1.5em}\textbf{Output2:} $\hat{\Pi}=\Pi_{k}$,\ $\hat{\bar{K}}=(R+\Lambda_{k}^{2})^{-1}\mathcal{\bar{K}}_{k+1}-\hat{K}$
    %\Statex \hspace{3.5em}Choose a stabilizer $\hat{K}+\bar{K}_{0}$ of system $(A+G,B)$.
    \Statex \hspace{3.5em}\textbf{while} $\| \bar{K}_{k+1}-\bar{K}_{k} \| > \epsilon$ \textbf{do}
    \Statex \hspace{3.5em}\qquad Solve $\{\Pi_{k},\mathcal{\bar{K}}_{k+1},\Lambda_{k}^{2}\}$ by (\ref{eq_Pi_K_lam})
    \Statex \hspace{3.5em}\qquad $k=k+1$
    \Statex \hspace{3.5em}\textbf{end while}
    \Statex \hspace{-1.5em}\textbf{Result:} \hspace{1em}Apply $u_{ik}\!=\!-\hat{K}x_{ik}\!-\!\hat{\bar{K}}x_{k}^{(N)}$, where $i\!=\!1,2,...,N$.
    \end{algorithmic}
\end{algorithm}
The notation $\hat{}$ in Algorithm \ref{alg_01} is used to indicate the estimated values of matrix coefficients and parameters, distinguishing them from the true values.
The Figure \ref{fig_logic} illustrates the logical diagram of the Algorithm \ref{alg_01}, providing a more intuitive understanding of the relationship between the different steps.
The effectiveness of the algorithm can be supported by theoretical guarantees provided by Theorems \ref{thm03} and \ref{thm04}.

\section{Simulation}
In this section, a numerical simulation is carried out to validate the effectiveness of the proposed algorithm. 
The large-scale population involves 500 agents, where the coefficients of each agent's dynamics are as follows
\begin{align*}
	&A=
	\begin{bmatrix}
		0.08    &  0.63\\
		0.39  &  0.26	
	\end{bmatrix} ,\
	B=
	\begin{bmatrix}
		0.10 \\
		0.16	
	\end{bmatrix} ,\\
	&G=
	\begin{bmatrix}
		0.10    &  0.05\\
		0.07  &  0.06	
	\end{bmatrix} ,\
	D=
	\begin{bmatrix}
		0.12    &  0.05\\
		0.11  &  0.12	
	\end{bmatrix} ,
\end{align*}
with $x_{ik} \in \mathbb{R}^{2}$, $u_{ik} \in \mathbb{R}$, and $w_{ik}\sim \mathcal{N}(0, 0.01)$. The initial state $x_{i0}$ is uniformly distributed on $[-6,0]\times[0,12] \subset \mathbb{R}^{2}$ with $\mathbb{E}[x_{i0}]=[6,-3]^{T}$. The parameters of the cost function (\ref{Ji_u}) are
\begin{align*}
	Q=
	\begin{bmatrix}
		2.00    &  -1.54\\
		-1.54    &  -0.12	
	\end{bmatrix}, \
	\Gamma=
	\begin{bmatrix}
		0.62    &  0.84\\
		0.80    &   0.54	
	\end{bmatrix} ,\
	R = -1.74,   
\end{align*}
where the eigenvalues of $Q$ are approximately $\lambda_{1}=2.8642$, $\lambda_{2}=-0.7442$. In this simulation, to implement Algorithm \ref{alg_01}, the discount factor is chosen as $\gamma=0.9$ and the control inputs of $\mathcal{A}_{1}$ and $\mathcal{A}_{2}$ are designed as
\begin{align}\label{eq_input}
	&K_{0}=
	\begin{bmatrix}
		0.05    &  -0.91
	\end{bmatrix} ,\
	\bar{K}_{0}=
	\begin{bmatrix}
		2.87    &  0.83
	\end{bmatrix} ,\cr
	&\xi_{ik}=\sum_{j=1}^{100}sin(w_{ik}^{j}),
\end{align}
where the frequencies $w_{ik}^{j}, i, j=1,\cdots,100$, are randomly selected from $[-100,100]$, and convergence criterion $\epsilon = 10^{-4}$. 

\begin{figure}[htbp]
    \centering
    \includegraphics[width=0.88\linewidth]{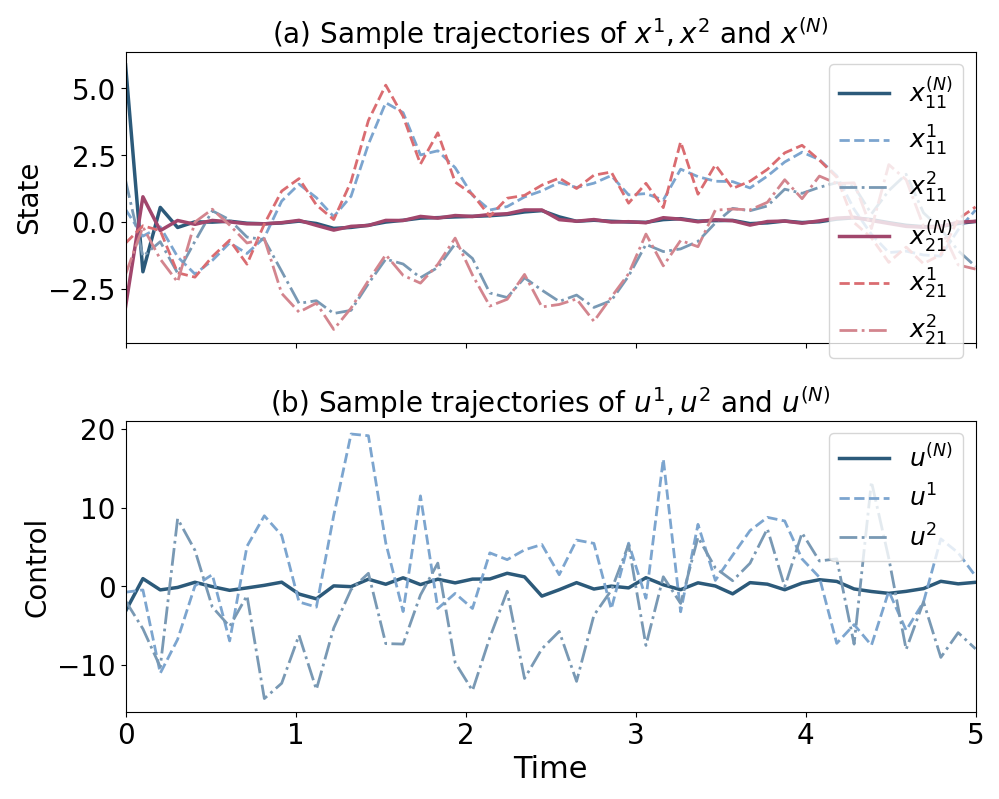}
    \caption{Real-time data collected from agent 1 and 2}
    \label{fig05}
\end{figure}

We apply the control input (\ref{eq_input}) to $\mathcal{A}_{1}$ and $\mathcal{A}_{2}$, and collect the dataset of the MF term as well as agents 1 and 2 after 50 iterations under the rank conditions (\ref{A5}) and (\ref{A6}). Fig. \ref{fig05} illustrates the states and control trajectories of the MF term and agents 1 and 2 under the effect of control input (\ref{eq_input}). 

\begin{figure}[htbp]
    \centering
    \includegraphics[width=0.88\linewidth]{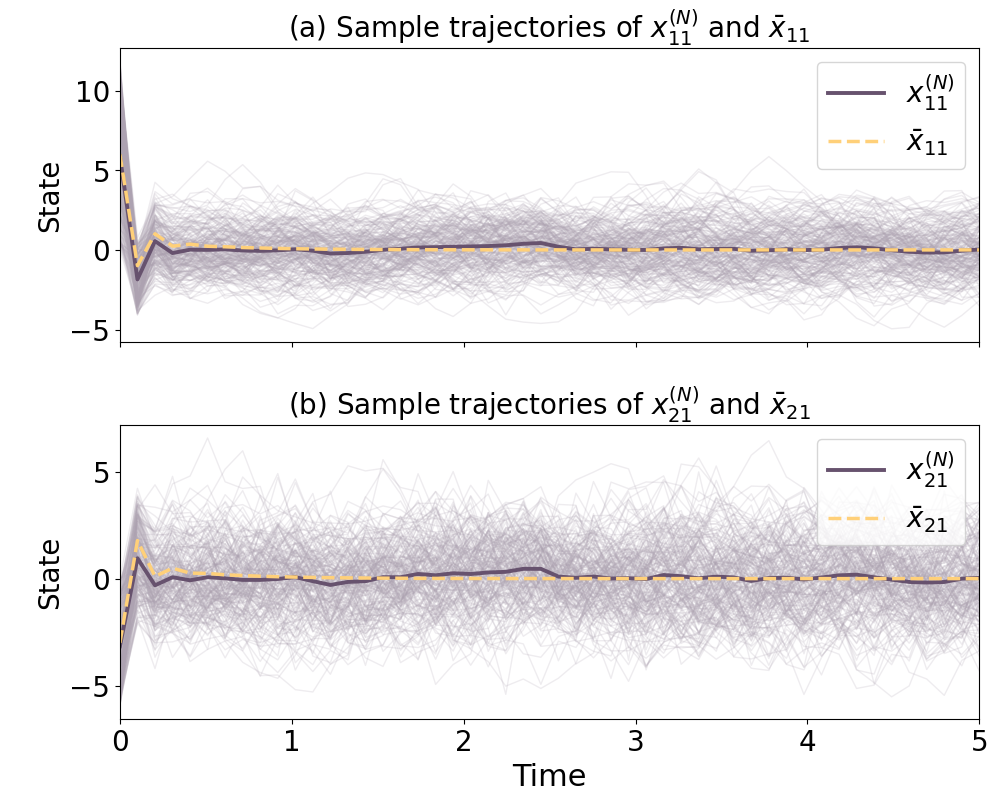}
    \caption{MF state trajectory comparison}
    \label{fig06}
\end{figure}

In addition, Fig. \ref{fig06} compares the MF state trajectory reconstructed from the learned gains with that obtained from the collected data, which shows that the learned controller can approximate the population behavior of the large-scale population system.

\begin{figure}[htbp]
	\centering
	\includegraphics[scale=0.3]{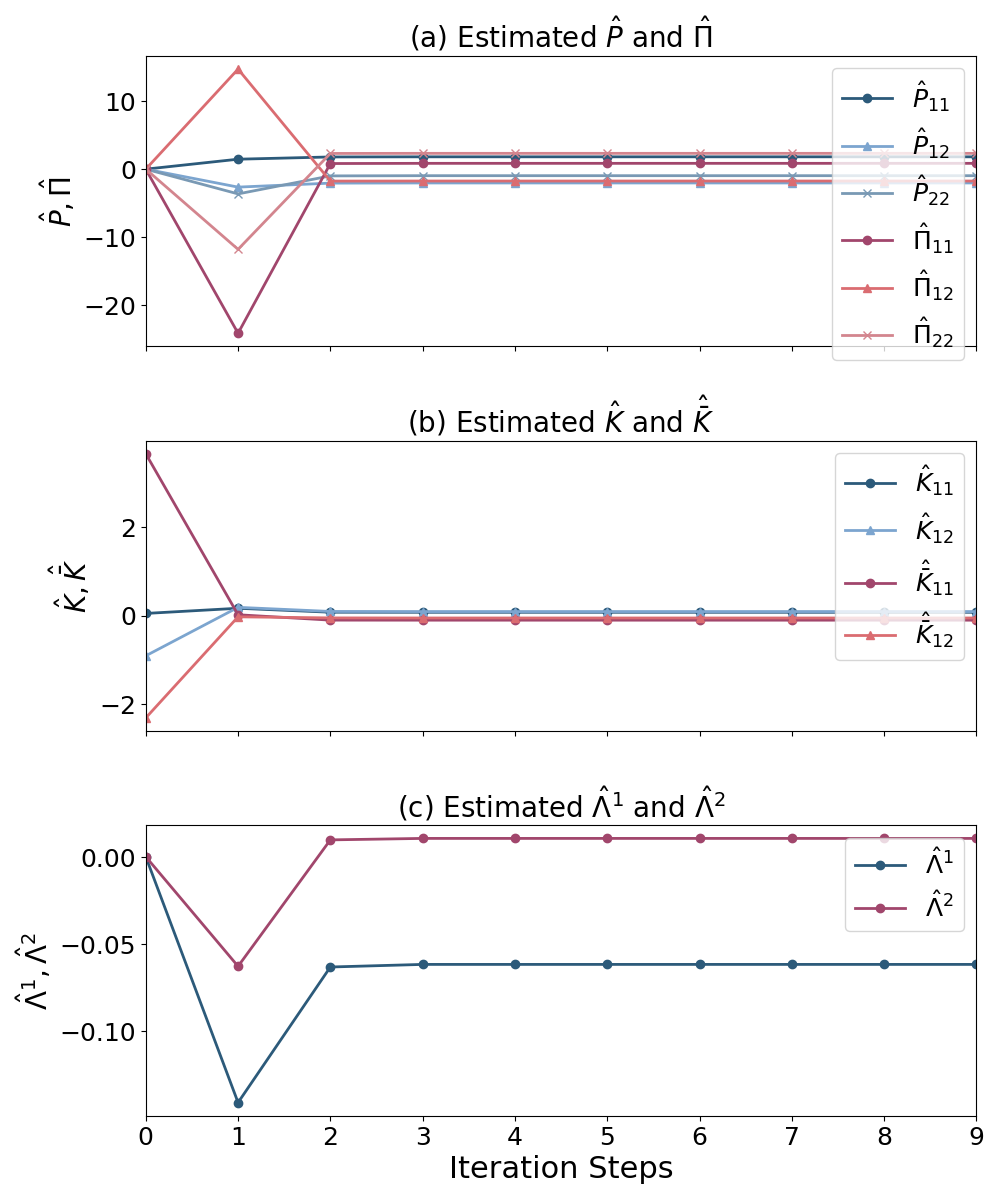}
    \caption{$\{\hat{P}, \hat{K},\hat{\Lambda}^{1}\}$ and $\{\hat{\Pi}, \hat{\bar{K}},\hat{\Lambda}^{2}\}$ of Algorithm \ref{alg_01}}
    \label{fig07}
\end{figure}

The convergence sequences $\{\hat{P}, \hat{K},\hat{\Lambda}^{1}\}$ and $\{\hat{\Pi}, \hat{\bar{K}},\hat{\Lambda}^{2}\}$ are shown in Fig. \ref{fig07}. Simulation results indicate that $\{\hat{P}, \hat{K},\hat{\Lambda}^{1}\}$ converges at the 3rd iteration under the convergence criterion $\epsilon$, while $\{\hat{\Pi}, \hat{\bar{K}},\hat{\Lambda}^{2}\}$ reaches convergence at the 3rd iteration. %Therefore, under certain conditions，even though the weight matrices are not positive semidefinite, the algorithm remains valid. 
Therefore, the algorithm remains valid and achieves convergence under certain conditions, despite the weight matrices not being positive semidefinite.

Tables \ref{tab01} and \ref{tab02} provide a summary of the final estimated parameters along with their corresponding relative errors. The model-based algorithm presented in Section~\ref{sec:Model-based-MF-social-control-design} was used to obtain the reference values $(\mathcal{P}^{*},\mathcal{K}^{*})$ and $(\Pi^{*},\bar{\mathcal{K}}^{*})$, which serve as the ground truth for comparison. The relative errors presented in the tables reflect the accuracy of the estimated values, indicating that the algorithm performs well in approximating the real values, with the estimates closely matching the theoretical results.

\begin{table}[htbp]
\renewcommand{\arraystretch}{1.2}
	\centering
	\footnotesize
	\caption{Estimates of  $\{\hat{P}, \hat{K}\}$}
	\label{tab01}
        \tabcolsep 5pt
	\begin{tabular}{cccccc}
		\toprule
		Model-free & Value & Model-based & Value & Error & Value \\
		\midrule
		$[\hat{P}_{11}]$ & 1.7614 & $[\mathcal{P}^{*}_{11}]$ &  1.7835 &      &        \\
		$[\hat{P}_{12}]$ & -2.0053 & $[\mathcal{P}^{*}_{12}]$ &  -2.0194 & $\dfrac{\|\mathcal{P}^{*}-\hat{P}\|_{2}}{\|\mathcal{P}^{*} \|_{2}}$ & 0.0112\\
		$[\hat{P}_{22}]$ & -0.9746& $[\mathcal{P}^{*}_{22}]$ &  -0.9653 &      &        \\
		\midrule
		$[\hat{K}_{11}]$ & 0.0723 & $[\mathcal{K}^{*}_{11}]$ & 0.0727 & \multirow{2}{*}{$\dfrac{\|\mathcal{K}^{*}-\hat{K}\|_{2}}{\|\mathcal{K}^{*}\|_{2}}$} & \multirow{2}{*}{0.0064}\\
		$[\hat{K}_{12}]$ & 0.0901 & $[\mathcal{K}^{*}_{12}]$ & 0.0905 &      &        \\
%		\midrule
%		$[\hat{\Lambda}^{1}]$ & -0.0617 & $[\Lambda^{*}_{11}]$ & 0.0727 & $\dfrac{\|\Lambda^{*}-\hat{\Lambda}\|_{2}}{\|\Lambda^{*}\|_{2}}$ & 0.0002\\
		\bottomrule
	\end{tabular}
\end{table}

\begin{table}[htbp]
\renewcommand{\arraystretch}{1.5}
	\centering
	\footnotesize
	\caption{Estimates of  $\{\hat{\Pi}, \hat{\bar{K}}\}$}
	\label{tab02}
        \tabcolsep 5pt
	\begin{tabular}{cccccc}
		\toprule
		Model-free & Value & Model-based & Value & Error & Value \\
		\midrule
		$[\hat{\Pi}_{11}]$ & 0.8556 & $[\Pi^{*}_{11}]$ &  0.8106 &      &        \\
		$[\hat{\Pi}_{12}]$ & -1.7479 & $[\Pi^{*}_{12}]$ &  -1.7897 & $\dfrac{\|\Pi^{*}-\hat{\Pi}\|_{2}}{\|\Pi^{*} \|_{2}}$ & 0.0278\\
		$[\hat{\Pi}_{22}]$ & 2.3274 & $[\Pi^{*}_{22}]$ &  2.2639 &      &        \\
		\midrule
		$[\hat{\bar{K}}_{11}]$ & -0.0978 & $[\mathcal{\bar{K}}^{*}_{11}]$ & 0.0976 & \multirow{2}{*}{$\dfrac{\|\mathcal{\bar{K}}^{*}-\hat{\bar{K}}\|_{2}}{\|\mathcal{\bar{K}}^{*}\|_{2}}$} & \multirow{2}{*}{0.0838}\\
		$[\hat{\bar{K}}_{12}]$ & -0.0249 & $[\mathcal{\bar{K}}^{*}_{12}]$ & -0.0508 &      &        \\
%		\midrule
%		$[\hat{\Lambda}^{1}]$ & -0.0617 & $[\Lambda^{*}_{11}]$ & 0.0727 & $\dfrac{\|\Lambda^{*}-\hat{\Lambda}\|_{2}}{\|\Lambda^{*}\|_{2}}$ & 0.0002\\
		\bottomrule
	\end{tabular}
\end{table}

\begin{figure}[htbp]
		\centering
		\includegraphics[scale=0.25]{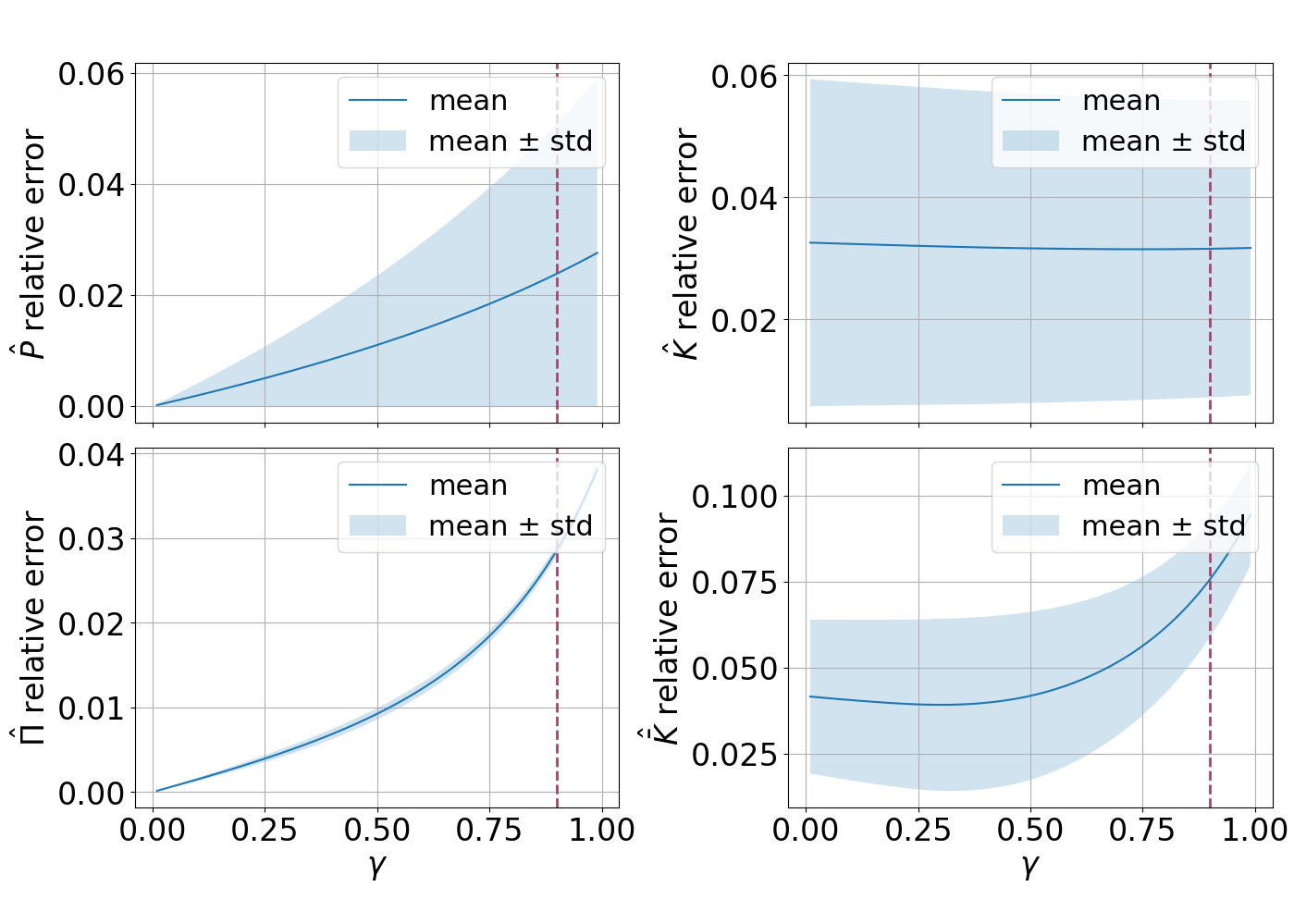}
	\caption{Relative estimation errors versus $\gamma$}
	\label{fig09}
\end{figure}

Furthermore, we investigate the influence of the discount factor $\gamma$. 
With all other parameters fixed, we vary $\gamma \in (0,1)$ and perform 200 runs for each value to compute the mean and standard deviation of the relative estimation errors, as shown in Fig. \ref{fig09}. 
It is observed that increasing $\gamma$ leads to a noticeable growth in the relative errors of $P$ and $\Pi$, and larger variances in the errors of $K$ and $\bar{K}$. 
This suggests that the learned feedback gains are more susceptible to randomness under finite data and noise.
Since smaller relative errors do not necessarily imply better long-term performance, we finally set $\gamma = 0.9$ as a compromise between estimation accuracy and closed-loop behavior.

\section{Conclusion}
In this paper, a data-driven RL algorithm is proposed to solve the decentralized asymptotic optimal control for discrete-time LQG-MF social control, where state coupling and the state and control weighting matrices are allowed not to be positive
semidefinite. 
Subsequently, a data-driven iterative optimization equation is derived through a system transformation method, which eliminates the dependence of the AREs on system dynamics. 
Simulation results validate the effectiveness of the proposed algorithm.
For future work, a possible extension is to explore leader-follower MF game and control problems.

\appendices
\section{Proofs of Lemma \ref{lem02} and Theorem \ref{thm01}} \label{app:A}
\begin{proof}[Proof of Lemma \ref{lem02}]
Under the assumption \ref{A1}, \ref{A2} and \ref{A4}, it can be shown that the Riccati equations (\ref{eq_P})-(\ref{eq_Pi}) admit a pair of $\gamma$-stabilizing solutions \(P,\Pi\) and a pair of negative definite solutions \(P^{-},\Pi^{-}\)(\cite{ran1993linear}). We have
\begin{align}
    \nonumber
    &J_{soc}(u)\\
    \nonumber
    =&\sum_{i=1}^{N}\mathbb{E}\Big\{\sum_{k=0}^{\infty}\gamma^{k}\big[\| x_{ik}-x_{k}^{(N)}\|_{Q}^{2}+\|x_{k}^{(N)}\|_{Q+Q_{\Gamma}}^{2}\\
    \nonumber
    &\!+\| u_{ik}-u_{k}^{(N)}\|_{R}^{2}+\|u_{k}^{(N)}\|_{R}^{2}\big]\Big\}\\
    \nonumber
    =&\sum_{i=1}^{N}\mathbb{E}\Big\{\|x_{i0}-x_{0}^{(N)}\|_{P^{-}}^{2}+\|x_{0}^{(N)}\|_{\Pi^{-}}^{2}\Big\}\\
    \nonumber
    &\!-\!\lim_{k \to \infty}\sum_{i=1}^{N}\mathbb{E}\Big\{\gamma^{k\!+\!1}\big[\|x_{i(k\!+\!1)}\!-\!x_{k\!+\!1}^{(N)}\|_{P^{\!-}}^{2}\!+\!\|x_{k\!+\!1}^{(N)}\|_{\Pi^{\!-}}^{2}\big]\Big\}\\
    \nonumber
    &\!+\!\sum_{i=1}^{N}\mathbb{E}\Big\{\sum_{k=0}^{\infty} \gamma^{k} \big[ \|u_{ik}\!-\!u_{k}^{(\!N\!)}\!+\!K(x_{ik}\!-\!x_{k}^{(\!N\!)})\|_{R+\gamma B^{T}\!\!P\!B}^{2}\\
    &\!+\|u_{k}^{(N)}+\bar{K}x_{k}^{(N)}\|_{R+\gamma B^{T}\!\Pi \!B}^{2}\big]\Big\}+q_{\infty}.
\end{align}    
Consider the following control set
\begin{align}
    \mathcal{U}_{c}^{'}=&\Big\{(u_{1},\dots,u_{N})|u_{ik}\in \mathbb{R}^{m} \ is\ adapted\ to\ \mathcal{F}_{i(k-1)},\cr
    &\mathbb{E}\{\sum_{k=0}^{\infty} \|x_{ik}\|^{2}\}< \infty,i=1,\dots,N\Big\}.
\end{align}
For any $u\in \mathcal{U}_{c}^{'}$ satisfying $J_{soc}(u)\leq NC$, we have
\begin{align}\label{eq_J_soc_u}
    \nonumber
    &J_{soc}(u)\\
    \nonumber
    =&\sum_{i=1}^{N}\mathbb{E}\Big\{\|x_{i0}-x_{0}^{(N)}\|_{P}^{2}+\|x_{0}^{(N)}\|_{\Pi}^{2}\Big\}\\
    \nonumber
    &\!+\!\sum_{i=1}^{N}\mathbb{E}\Big\{\sum_{k=0}^{\infty} \gamma^{k} \big[ \|u_{ik}\!-\!u_{k}^{\!(\!N\!)}\!+\!K(x_{ik}\!-\!x_{k}^{\!(\!N\!)})\|_{R+\gamma B^{T}\!\!P\!B}^{2}\\
    &\!+\|u_{k}^{(N)}+\bar{K}x_{k}^{(N)}\|_{R+\gamma B^{T}\!\Pi \!B}^{2}\big]\Big\}+q_{\infty} \leq NC.
\end{align}
Denote $v_{k}^{(N)} = u_{k}^{(N)}+\bar{K}x_{k}^{(N)}$.
From (\ref{xi_k+1}),
\begin{align}
    x_{k+1}^{(N)}\!=\!(A\!+\!G\!-\!B\bar{K})x_{k}^{\!(N)}\!+\!Bv_{k}^{\!(N)}\!+\!\frac{1}{N}\!\sum_{i=1}^{N}\!Dw_{ik}.
\end{align}
By Huang \cite{huang2010large}, there exist constants $C_{1},C_{2}>0$ such that
\begin{align}
    \mathbb{E}\Big\{ \sum_{k=0}^{\infty}\gamma^{k}\|x_{k}^{(N)}\|^{2}\Big\} \leq C_{1}\mathbb{E}\Big\{ \sum_{k=0}^{\infty}\gamma^{k}\|v_{k}^{(N)}\|^{2}\Big\}+C_{2}.
\end{align}
Combining this inequality with (\ref{eq_J_soc_u}) yields
\begin{align}\label{eq_x_N_u_N}
    \nonumber
    &\sum_{i=1}^{N}\mathbb{E}\Big\{ \sum_{k=0}^{\infty}\gamma^{k}\big[\|x_{k}^{(N)}\|^{2}+\|u_{k}^{(N)}\|^{2}\big]\Big\}\\
    \nonumber
    =&N\mathbb{E}\Big\{ \sum_{k=0}^{\infty}\gamma^{k}\big[\|x_{k}^{(N)}\|^{2}+\|v_{k}^{(N)}-\bar{K}x_{k}^{(N)}\|^{2}\big]\Big\}\\
    \leq& NC_{3}\mathbb{E}\Big\{ \sum_{k=0}^{\infty}\gamma^{k}\|v_{k}^{(N)}\|^{2}\Big\}+NC_{4} \leq NC.
\end{align}
Similarly, we obtain
\begin{align}
    \sum_{i=1}^{N}\mathbb{E}\Big\{\sum_{k=0}^{\infty}\gamma^{k}\big[\| x_{ik}\!-\!x_{k}^{\!(N)}\|^{2}\!+\!\| u_{ik}\!-\!u_{k}^{\!(N)}\|^{2}\big]\Big\}\!\leq NC.
\end{align}
Consequently, from the above and (\ref{eq_x_N_u_N}), it follows that
\begin{align}
    \sum_{i=1}^{N}\mathbb{E}\Big\{\sum_{k=0}^{\infty}\gamma^{k}\big[\| x_{ik}\|^{2}+\| u_{ik}\|^{2}\big]\Big\}\leq NC.
\end{align}
The remaining part of the proof is similar to that of Lemma \ref{lem01}.  
For the case that initial state $\{x_{i0}\}$ share the same variance, then from (\ref{eq_J_soc_hat_u}), the asymptotic average social optimum is given by $\mathbb{E}\big[\|x_{i0}-\bar{x}_{0}\|_{P}^{2}+\|\bar{x}_{0}\|_{\Pi}^{2}\big]+q_{\infty}$.
\end{proof}

\begin{proof}[Proof of Theorem \ref{thm01}]
Assumptions \ref{A3} and \ref{A4} specify different properties of the matrix $Q$. This leads to a different proof for part a), while the proofs for parts b) and c) are similar. Specifically:\\
a) We divide the proof of part a) into two parts:\\
\textbf{(i) The case that $Q \ge 0$ and $R > 0$} (under the assumption \ref{A3}).
	To simplify the proof, we abbreviate equation (\ref{eq_P_k}) as follows
	\begin{align}\label{eq_A_S}
		P_{k}=\gamma A_{k}^{T}P_{k}A_{k}+S_{k}, \ where \ S_{k}=K_{k}^{T}RK_{k}+Q.
	\end{align}
Given an arbitrary $\gamma$-stabilizing feedback gain matrix $K_{0}$, assume that $K_{k}$ $\gamma$-stabilizes $(A, B)$ for $k \geq 1$. We can prove by contradiction that $(\sqrt{\gamma}A_{k},\sqrt{S_{k}})$ is detectable.
According to Lemma \ref{lem03}, assume that there exists a non-zero symmetric matrix $\mathcal{X}$ such that
\begin{align}\label{eq_A_k_S_k}
	\left\{
	\begin{aligned}
		&\mathcal{X} - \gamma A_{k}^{T}\mathcal{X}A_{k}=\lambda \mathcal{X}, \ \lvert \lambda \rvert \geq 1,\\
		&\sqrt{S_{k}}\mathcal{X} = \mathbf{0}.   
	\end{aligned}
	\right.    
\end{align}
Due to $R > 0$, we have $K_{k}\mathcal{X}=\mathbf{0}$ and $\sqrt{Q}\mathcal{X}=\mathbf{0}$. Equation (\ref{eq_A_k_S_k}) can be simplified as follows
\begin{align}\label{eq_A_Q}
	\left\{
	\begin{aligned}
		&\mathcal{X} - \gamma A^{T}\mathcal{X}A=\lambda \mathcal{X}, \ \lvert \lambda \rvert \geq 1,\\
		&\sqrt{Q}\mathcal{X} = \mathbf{0},
	\end{aligned}
	\right.    
\end{align}
which implies that $(\sqrt{\gamma}A,\sqrt{Q})$ is not detectable, contradicting Assumption~\ref{A3}. Therefore, $(A_{k},\sqrt{S_{k}})$ is detectable.
Additionally, $K_{k}$ is a $\gamma$-stabilizer of $(A,B)$. According to Theorem 3 in \cite{2019zhang_dt_mfc}, $P_{k} \geq 0$ is the unique positive semidefinite solution to equation (\ref{eq_P_k}).
In order to demonstrate that $K_{k+1}$ serves as a $\gamma$-stabilizer for $(A,B)$, we rewrite equation (\ref{eq_A_S}) as follows
\begin{align}
	P_{k}=\gamma A_{k+1}^{T}P_{k}A_{k+1}+\Tilde{S}_{k+1},     
\end{align}
where $\Tilde{S}_{k+1} = (K_{k}-K_{k+1})^{T}(R+\gamma B^{T}P_{k}B)(K_{k}-K_{k+1})+K_{k+1}^{T}RK_{k+1}+Q$.
Based on the above derivation, combined with the positive semidefinite condition of $P_{k}\geq 0$ and the conclusion of Theorem 3 in reference \cite{2019zhang_dt_mfc}, it can be proven that $(\sqrt{\gamma}A_{k+1},\sqrt{\Tilde{S}_{k+1}})$ is detectable, and $K_{k+1}$ constitutes a $\gamma$-stabilizer for the system $(A,B)$.
Hence, $K_{k}$ is a $\gamma$-stabilizer of the system $(A,B)$, which implies that $\sqrt{\gamma}A_{k}$ is Schur.\\
\textbf{(ii) The case that $Q$ and $R$ are symmetric matrices} (under the assumption \ref{A4}). For $k=0$, the matrix $\sqrt{\gamma}A_{0}$ is Schur due to the $\gamma$-stabilizing feedback gain matrix $K_{0}$. Thus, equation (\ref{eq_P_k}) transforms to
%当k=0时，由于K_{0}是系统(A, B)的稳定器，所以A_{0}是舒尔的，公式(10)可以写为
\begin{align}\label{eq_P_0_A_0} 
P_{0}=\gamma A_{0}^{T}P_{0}A_{0}+K_{0}^{T}RK_{0}+Q,
\end{align}
For $k \geq 1$, assuming that $\sqrt{\gamma}A_{k}$ is Schur, equation (\ref{eq_P_0_A_0}) can be rewritten as
%当k大于等于1时，假定A_{k}是舒尔的，公式\ref{eq_P_k_0}可以写为
\begin{align}\label{eq_P_0_A_k} 
\nonumber
P_{0}=&\gamma A_{k}^{T}P_{0}A_{k}-(K_{k}-K_{0})^{T}(R+\gamma B^{T}P_{0}B)(K_{k}-K_{0})\\
&+K_{k}^{T}RK_{k}+Q,
\end{align}
then we have
\begin{align}\label{P_k-P_0}
	P_{k}-P_{0}=&\sum_{n=0}^{\infty}\gamma^{n} (A_{k}^{T})^{n}(K_{k}-K_{0})^{T}(R+\gamma B^{T}P_{0}B)\cr
    &\times(K_{k}-K_{0})(A_{k})^{n} \geq 0.
\end{align} 
Next, by contradiction, we show the matrix $\sqrt{\gamma}A_{k+1}$ is Schur and rewrite equation (\ref{P_k-P_0}) as
 %之后，我们通过反证法证明A_{k+1}是舒尔的，改写公式(\ref{P_k-P_0)
\begin{align}\label{P_k_A_k+1}
	P_{k}&-P_{0}=\gamma A_{k+1}^{T}(P_{k}-P_{0})A_{k+1}\cr
    &+(K_{k+1}-K_{0})^{T}(R+\gamma B^{T}P_{0}B)(K_{k+1}-K_{0})\cr
    &+(K_{k}-K_{k+1})^{T}(R+\gamma B^{T}P_{k}B)(K_{k}-K_{k+1}).
\end{align}
Assume $\sqrt{\gamma}A_{k+1}\mathtt{z}=\lambda_{i} \mathtt{z}$, for $|\lambda_{i}| \geq 1$ and $\mathtt{z} \neq 0$, then we have
\begin{align}\label{zA_k+1}
	&\gamma \mathtt{z}^{T}\!A_{k+1}^{T}\!(P_{k}\!-\!P_{0})A_{k+1}\mathtt{z}\!-\!\mathtt{z}^{T}\!(P_{k}\!-\!P_{0})\mathtt{z}\cr
    \!\!\!\!= &-\!\mathtt{z}^{T}\!(K_{k+1}\!-\!K_{0})^{T}\!(R\!+\!\gamma B^{T}\!P_{0}B)(K_{k+1}\!-\!K_{0})\mathtt{z}\cr
    &-\!\mathtt{z}^{T}\!(K_{k}\!-\!K_{k+1})^{T}\!(R\!+\!\gamma B^{T}\!P_{k}B)(K_{k}\!-\!K_{k+1})\mathtt{z}\!\leq\!0.
\end{align}    	
Substituting $\sqrt{\gamma}A_{k+1}\mathtt{z}=\lambda_{i} \mathtt{z}$ into the equation (\ref{zA_k+1}), we have
\begin{align}\label{zA_k+1_lam}
    &\gamma \mathtt{z}^{T}A_{k+1}^{T}(P_{k}-P_{0})A_{k+1}\mathtt{z}-\mathtt{z}^{T}(P_{k}-P_{0})\mathtt{z}\cr
    =&(\lambda_{i}^{2}-1)\mathtt{z}^{T}(P_{k}-P_{0})\mathtt{z} \geq 0.
\end{align}
Thus, combining inequality (\ref{zA_k+1_lam}) with (\ref{zA_k+1}), we conclude that
\begin{align}
	-&\mathtt{z}^{T}\!(K_{k}\!-\!K_{k\!+\!1})^{T}\!(R\!+\!\gamma B^{T}\!P_{k}B)(K_{k}\!-\!K_{k\!+\!1})\mathtt{z}\cr
    &-\!\mathtt{z}^{T}\!(K_{k\!+\!1}\!-\!K_{0})^{T}\!(R\!+\!\gamma B^{T}\!P_{0}B)(K_{k\!+\!1}\!-\!K_{0})\mathtt{z}\!=\!0,
\end{align}	   
which gives rise to $(K_{k}-K_{k+1})\mathtt{z}=0$. Consequently, we can get $\sqrt{\gamma}A_{k+1}\mathtt{z}=\sqrt{\gamma}A_{k}\mathtt{z}=\lambda_{i} \mathtt{z}$. It contradicts the induction assumption.
Therefore, by mathematical induction, we have ultimately proven that $\sqrt{\gamma}A_{k}$ is Schur.
%所以通过数学归纳法，我们最终证明了A_{k}是舒尔的。

b) Using the result in part a), $\sqrt{\gamma}A_{k}$ is Schur. We rewrite equation (\ref{eq_P}) as
\begin{align}\label{P_A_k}
    \nonumber
	P=& \gamma A_{k}^{T}PA_{k}-(K_{k}-K)^{T}(R+\gamma B^{T}PB)(K_{k}-K)\\ 
    &+K_{k}^{T}RK_{k}+Q.
\end{align}
Then we have
\begin{align}\label{P_k-P}
	&P_{k}-P\cr
	\!=&\gamma A_{k}^{T}(P_{k}\!-\!P)A_{k}\!+\!(K_{k}\!-\!K)^{T}(R\!+\!\gamma B^{T}PB)(K_{k}\!-\!K)\cr
	\!=&\!\sum_{n=0}^{\infty} \!\gamma^{n}(A_{k}^{T})^{n}\!(K_{k}\!-\!K)^{T}\!(R\!+\!\gamma B^{T}\!PB)\!(K_{k}\!-\!K)\!(A_{k})^{n}\!,
\end{align} 
which yields $P_{k} \geq P$, as $R+\gamma B^{T}PB > 0$.
By equation (\ref{eq_P_k})
\begin{align}
	P_{k+1}&=\gamma A_{k+1}^{T}P_{k+1}A_{k+1}+K_{k+1}^{T}RK_{k+1}+Q,\\
    \nonumber
	P_{k}=&\gamma A_{k+1}^{T}P_{k}A_{k+1}+(K_{k}-K_{k+1})^{T}(R+\gamma B^{T}P_{k}B)\\
    &\times(K_{k}-K_{k+1})+K_{k+1}^{T}RK_{k+1}+Q.
\end{align}
Then we can get
\begin{align}
\nonumber
P_{k}\!-\!P_{k+1}
= &\sum_{n=0}^{\infty}\gamma^{n}(A_{k+1}^{T})^{n}(K_{k}\!-\!K_{k+1})^{T}\\
&\times\!(R\!+\!\gamma B^{T}\!P_{k}B)(K_{k}\!-\!K_{k+1})(A_{k+1})^{n},
\end{align}    
which yields $P_{k} \geq P_{k+1}$, as $R+\gamma B^{T}P_{k}B > 0$. Combining with the obtained result, we have $P_{k} \geq P_{k+1} \geq P$.

c) Since $K_{k+1}$ is uniquely determined by \eqref{eq_K_k}, the convergence of the sequence $\{P_{k}\}_{0}^{\infty}$ implies that $\{K_{k}\}_{0}^{\infty}$ also converges. It follows from b) that $\{P_{k}\}_{0}^{\infty}$ is monotonically decreasing sequence and has a lower bound $P$, leading to $\lim_{k \rightarrow \infty}P_{k}=P$. Hence, the proof is complete.
\end{proof}

\section{Proofs of Theorem \ref{thm02} and Theorem \ref{thm03}} \label{app:B}

\begin{proof}[Proof of Theorem \ref{thm02}]
Assumptions \ref{A3} and \ref{A4} define different properties of $Q$, resulting in a distinct proof for part a), with parts b) and c) following similarly. \\
a) We divide the proof of part a) into two parts:\\
(i) \textbf{The case that $Q \geq 0$ and $R > 0$} (under the assumption \ref{A3}). To simplify the proof, we abbreviate equation (\ref{eq_Pi_k}) as follows
	\begin{align}\label{eq_A_S_bar}
		\Pi_{k}=\gamma \bar{A}_{k}^{T}\Pi_{k}\bar{A}_{k}+\bar{S}_{k}, 
	\end{align}
    where $\bar{S}_{k}=(K+\bar{K}_{k})^{T}R(K+\bar{K}_{k})+Q +Q_{\Gamma}$.
	Given that $\bar{K}_{0}+K$ be any $\gamma$-stabilizing feedback gain matrix, and assuming that $\bar{K}_{k}+K$ is a $\gamma$-stabilizer of $(A+G,B)$ for $k \geq 1$. We can prove by contradiction that $(\sqrt{\gamma}\bar{A}_{k},\sqrt{\bar{S}_{k}})$ is detectable.
	According to Theorem 3 in \cite{2019zhang_dt_mfc}, $\Pi_{k} \geq 0$ is the unique positive semidefinite solution to equation (\ref{eq_Pi_k}).
	
	In order to demonstrate that $\bar{K}_{k+1}+K$ serves as a $\gamma$-stabilizer for $(A+G,B)$, we rewrite equation (\ref{eq_A_S_bar}) as follows
	\begin{align}
		\Pi_{k}=\gamma \bar{A}_{k+1}^{T}\Pi_{k}\bar{A}_{k+1}+\Tilde{\bar{S}}_{k+1},    
	\end{align}
	where $\Tilde{\bar{S}}_{k+1} = (\bar{K}_{k}-\bar{K}_{k+1})^{T}(R+\gamma B^{T}\Pi_{k}B)(\bar{K}_{k}-\bar{K}_{k+1})+(\bar{K}_{k+1}+K)^{T}R(\bar{K}_{k+1}+K)+Q+Q_{\Gamma}$,
	and it can be proven that $(\sqrt{\gamma} \bar{A}_{k+1},\sqrt{\Tilde{\bar{S}}_{k+1}})$ is detectable, combined with $\Pi_{k} \geq 0$ and the conclusion of Theorem 3 in reference \cite{2019zhang_dt_mfc}, which results in $\bar{K}_{k+1}+K$ being a $\gamma$-stabilizer of $(A+G,B)$.	
	Hence, $\bar{K}_{k}+K$ is a $\gamma$-stabilizer of the system $(A+G,B)$, which implies that $\sqrt{\gamma}\bar{A}_{k}$ is Schur.

(ii) \textbf{The case that $Q$ and $R$ are symmetric matrices} (under the assumption \ref{A4}). For $k=0$, the matrix $\sqrt{\gamma}\bar{A}_{0}$ is Schur because $\bar K_{0}+K$ $\gamma$-stabilizes $(A+G,B)$.
For $k \geq 1$, assuming that $\sqrt{\gamma}\bar{A}_{k}$ is Schur, and combined with equation (\ref{Pi_A_G}), we have
\begin{align}\label{Pi_k-Pi_0}
\Pi_{k}\!-\!\Pi_{0}=&\!\sum_{n=0}^{\infty}\gamma^{n}(\bar{A}_{k}^{T})^{n}(\bar{K}_{k} \!-\! \bar{K}_{0})^{T}\cr
&\!\times(R \!+\! \gamma B^{T}\Pi_{0}B)(\bar{K}_{k} \!-\! \bar{K}_{0})(\bar{A}_{k})^{n} \geq  0,
\end{align} 
Next, we show that the matrix $\sqrt{\gamma}\bar{A}_{k+1}$ is Schur by contradiction. 
equation (\ref{Pi_k-Pi_0}) can be rewriten as
\begin{align}\label{Pi_k_A_k+1}
	\Pi_{k}&-\Pi_{0}=\gamma \bar{A}_{k+1}^{T}(\Pi_{k}-\Pi_{0})\bar{A}_{k+1}\cr
    &+(\bar{K}_{k+1}-\bar{K}_{0})^{T}(R+\gamma B^{T}\Pi_{0}B)(\bar{K}_{k+1}-\bar{K}_{0})\cr
    &+(\bar{K}_{k}-\bar{K}_{k+1})^{T}(R+\gamma B^{T}\Pi_{k}B)(\bar{K}_{k}-\bar{K}_{k+1}).
\end{align}
Assume $\sqrt{\gamma}\bar{A}_{k+1}\mathtt{z}=\bar{\lambda}_{i} \mathtt{z}$, for $|\bar{\lambda}_{i}| \geq 1$ and $\mathtt{z} \neq 0$, we have 
 %之后，我们通过反证法证明A_{k+1}是舒尔的，改写公式(\ref{P_k-P_0)
\begin{align}\label{zA_k+1_bar}
     &\gamma \mathtt{z}^{T}\!\bar{A}_{k+1}^{T}\!(\Pi_{k}\!-\!\Pi_{0})\bar{A}_{k+1}\mathtt{z}\!-\!\mathtt{z}^{T}\!(\Pi_{k}\!-\!\Pi_{0})\mathtt{z}\cr
    \!\!\!=&-\!\mathtt{z}^{T}\!(\bar{K}_{k\!+\!1}\!-\!\bar{K}_{0})^{T}\!(R\!+\!\gamma B^{T}\!\Pi_{0}B)(\bar{K}_{k\!+\!1}\!-\!\bar{K}_{0})\mathtt{z}\cr
    &-\!\mathtt{z}^{T}\!(\bar{K}_{k}\!-\!\bar{K}_{k\!+\!1})^{T}\!(R\!+\!\gamma B^{T}\!\Pi_{k}B)(\bar{K}_{k}\!-\!\bar{K}_{k\!+\!1})\mathtt{z}\!\leq\! 0,
\end{align}    	
Substituting $\sqrt{\gamma}\bar{A}_{k+1}\mathtt{z}=\bar{\lambda}_{i} \mathtt{z}$ into the equation (\ref{zA_k+1_bar}), we have
\begin{align}\label{zA_k+1_lam_bar}
    &\gamma \mathtt{z}^{T}\bar{A}_{k+1}^{T}(\Pi_{k}-\Pi_{0})\bar{A}_{k+1}\mathtt{z}-\mathtt{z}^{T}(\Pi_{k}-\Pi_{0})\mathtt{z}\cr
    =&(\bar{\lambda}_{i}^{2}-1)\mathtt{z}^{T}(\Pi_{k}-\Pi_{0})\mathtt{z} \geq 0,
\end{align}
Thus, combining (\ref{zA_k+1_lam_bar}) with (\ref{zA_k+1_bar}), we conclude that $(\bar{K}_{k}\!-\!\bar{K}_{k\!+\!1})\mathtt{z}\!=\!0$. Consequently, we can get $\sqrt{\gamma}\bar{A}_{k\!+\!1}\mathtt{z}\!=\!\sqrt{\gamma}\bar{A}_{k}\mathtt{z}\!=\!\bar{\lambda}_{i} \mathtt{z}$. It contradicts the induction assumption.
Therefore, by mathematical induction, we have ultimately proven that $\sqrt{\gamma} \bar{A}_{k}$ is Schur.

b) Using the result in part a), $\sqrt{\gamma}\bar{A}_{k}$ is Schur. We rewrite equation (\ref{eq_Pi}) and then we have
\begin{align}\label{Pi_k-Pi}
    \nonumber
	\Pi_{k}-\Pi=&\sum_{n=0}^{\infty} \gamma^{n}(\bar{A}_{k}^{T})^{n}(\bar{K}_{k}-\bar{K})^{T}\\
    &\times(R+\gamma B^{T}\Pi B)(\bar{K}_{k}-\bar{K})(\bar{A}_{k})^{n},
\end{align}   
We can get $\Pi_{k} \geq \Pi$ by $R+\gamma B^{T}\Pi B > 0$. By equation (\ref{eq_Pi_k}), we have
\begin{align}
    \nonumber
	\Pi_{k}\!-\!\Pi_{k\!+\!1}=&\sum_{n=0}^{\infty}\gamma^{n}(\bar{A}_{k+1}^{T})^{n}(\bar{K}_{k}\!-\!\bar{K}_{k+1})^{T}\\
    &\times\!(R\!+\!\gamma B^{T}\Pi_{k}B)(\bar{K}_{k}\!-\!\bar{K}_{k\!+\!1})(\bar{A}_{k\!+\!1})^{n},
\end{align}   	
which yields $\Pi_{k} \geq \Pi_{k+1}$, as $R+\gamma B^{T}\Pi_{k}B > 0$. Combining with the previously obtained result, we have $\Pi_{k} \geq \Pi_{k+1} \geq \Pi$.
	
c) Since $\bar{K}_{k+1}$ is the unique solution of equation (\ref{eq_K_k_bar}), proving the convergence of the sequence $\{\Pi_{k}\}_{0}^{\infty}$ would ensure that $\{\bar{K}_{k}\}_{0}^{\infty}$ also converge. It follows from b) that $\{\Pi_{k}\}_{0}^{\infty}$ is monotonically decreasing sequence and has a lower bound $\Pi$, leading to $\lim_{k \rightarrow \infty}\Pi_{k}=\Pi$. Hence, the proof is complete.
\end{proof}

\begin{proof}[Proof of Theorem \ref{thm03}]
To prove that equation (\ref{eq_P_K_lam}) has a unique solution, we need to show that the matrix $\mathfrak{A}_{k}^{1}$ is of column full rank. The convergence result follows Theorem \ref{thm01}. 
Next, we show that $\mathfrak{A}_{k}^{1}$ is of column full rank.

Assume that there exists a vector $\mathcal{S}=[svec(\mathcal{S}_{1}),vec(\mathcal{S}_{2}),\\svec(\mathcal{S}_{3})]^{T}$, such that
\begin{align}
    \mathfrak{A}_{k}^{1}\mathcal{S}=\boldsymbol{0},
\end{align}
where $\mathcal{S}_{1} \!\in\! \mathbb{R}^{\frac{n(n+1)}{2}},\mathcal{S}_{2} \!\in\! \mathbb{R}^{m\times n},\mathcal{S}_{3} \!\in\! \mathbb{R}^{\frac{m(m+1)}{2}}$. Then we have
\begin{align}\label{eq_prf_1}
    \nonumber
    &\!\!(\mathcal{I}_{\Delta x\Delta x}\!-\!\gamma \mathcal{I}_{\Delta x\Delta x}^{'})svec(\mathcal{S}_{1})\!+\!\mathcal{I}_{\Delta x\Delta x}svec(K_{k}^{T}\!\mathcal{S}_{2}\!+\!\mathcal{S}_{2}^{T}\!K_{k}\\
    &\!\!-\!K_{k}^{T}\!\mathcal{S}_{3}K_{k})\!+\!2\mathcal{I}_{\Delta x\Delta u}vec(\mathcal{S}_{2})\!+\!\mathcal{I}_{\Delta u\Delta u}svec(\mathcal{S}_{3})\!=\!\boldsymbol{0}.
\end{align}   
According to the equation (\ref{delta_K+1_k_1}), it gives
\begin{align}\label{eq_prf_2}
    \nonumber
    &\!\!(\gamma \mathcal{I}_{\Delta x\Delta x}^{'}\!-\!\mathcal{I}_{\Delta x\Delta x})svec(\mathcal{S}_{1})\!=\!\mathcal{I}_{\Delta x\Delta x}svec(\gamma A^{T}\!\mathcal{S}_{1}A\!-\!\mathcal{S}_{1})\\
    &\!\!+\!2\gamma \mathcal{I}_{\Delta x\Delta u}vec(B^{T}\!\mathcal{S}_{1}A)\!+\!\gamma \mathcal{I}_{\Delta u\Delta u}svec(B^{T}\!\mathcal{S}_{1}B).
\end{align} 
Combining equation (\ref{eq_prf_1}) and (\ref{eq_prf_2}), we can get
\begin{align}
  \!\!\!\!\mathcal{I}_{\Delta x\Delta x} svec(\mathcal{U}_{1})\!+\!\!2\mathcal{I}_{\Delta x\Delta u} vec(\mathcal{U}_{2})\!+\!\!\mathcal{I}_{\Delta u\Delta u} svec(\mathcal{U}_{3})\!=\!\boldsymbol{0},  
\end{align}
where
\begin{align*}
    \nonumber
    &\mathcal{U}_{1}=\gamma A^{T}\mathcal{S}_{1}A-\mathcal{S}_{1}-K_{k}^{T}\mathcal{S}_{2}-\mathcal{S}_{2}^{T}K_{k}+K_{k}^{T}\mathcal{S}_{3}K_{k},\\
    \nonumber
    &\mathcal{U}_{2}=\gamma B^{T}\mathcal{S}_{1}A-\mathcal{S}_{2},\\
    \nonumber
    &\mathcal{U}_{3}=\gamma B^{T}\mathcal{S}_{1}B-\mathcal{S}_{3}.
\end{align*}
Based on the rank condition (\ref{eq_rank_1}), we can derive that $\mathcal{U}_{1}=\mathcal{U}_{2}=\mathcal{U}_{3}=\boldsymbol{0}$. Then we have
\begin{align}
    \mathcal{S}_{1}-\gamma A_{k}^{T}\mathcal{S}_{1}A_{k}=\boldsymbol{0},
\end{align}
Since Theorem \ref{thm01} previously proved that $\sqrt{\gamma} A_{k}$ is Schur. According to \cite{hewer1971iterative}, we can conclude that $\mathcal{S}_{1}=0$, which in turn implies that $\mathcal{S}_{2}=\mathcal{S}_{3}=0$.
Thus, we have $\mathcal{S}=\boldsymbol{0}$,and so $\mathfrak{A}_{k}^{1}$ has full column rank.
\end{proof}

\vspace{-2\baselineskip}
\begin{IEEEbiography}[{\includegraphics[width=1in,height=1.25in,clip,keepaspectratio]{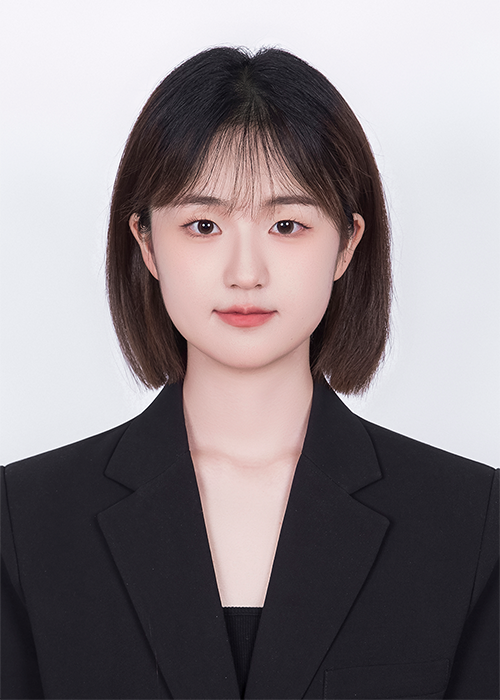}}]{Hanfang Zhang} received the B.E. degree in Automation from Qinghai University, Qinghai, China, in 2022.
She is currently pursuing the Ph.D. degree in Control Theory and Control Engineering from Shandong University, Jinan, China.

Her research interests include reinforcement learning, multi-agent systems, and mean field games.
\end{IEEEbiography}
\vspace{-2\baselineskip}
\begin{IEEEbiography}[{\includegraphics[width=1in,height=1.25in, clip,keepaspectratio]{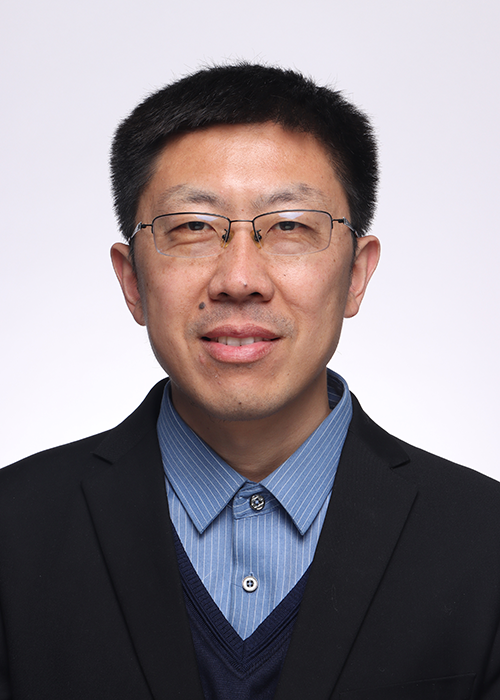}}]{Bing-Chang Wang} (SM’19) received the Ph.D. degree in System Theory from Academy of Mathematics and Systems Science, Chinese Academy of Sciences, Beijing, China, in 2011. From September 2011 to August 2012, he was with Department of Electrical and Computer Engineering, University of Alberta, Canada, as a Postdoctoral Fellow. From September 2012 to September 2013, he was with School of Electrical Engineering and Computer Science, University of Newcastle, Australia, as a Research Academic. From October 2013, he has been with School of Control Science and Engineering, Shandong University, China, and now is a
Professor. 

He held visiting appointments as a Research Associate with Carleton University, Canada, from November 2014 to May 2015, and with the Hong Kong Polytechnic University from November 2016 to January 2017. His current research interests include mean field games, stochastic control, multi-agent systems, and reinforcement learning.
\end{IEEEbiography}
\vspace{-2\baselineskip}
\begin{IEEEbiography}[{\includegraphics[width=1in,height=1.25in,clip,keepaspectratio]{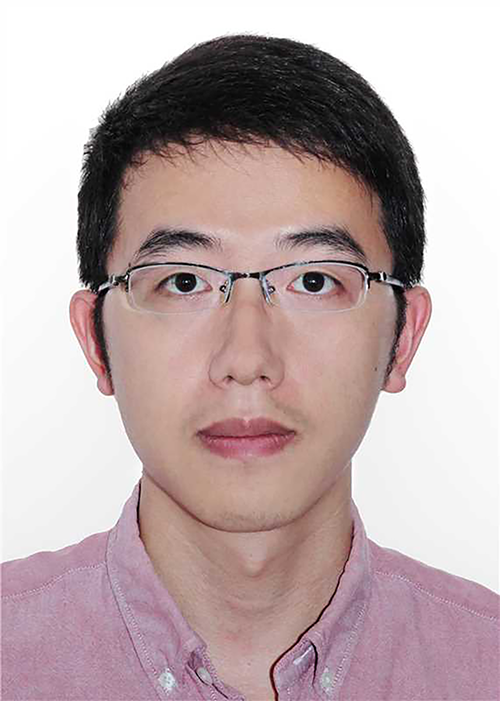}}]{Shuo Chen} received the Ph.D. degree in computer science from Nanyang Technological University, Singapore, in 2020, the M.Eng. degree in communication engineering from Nanyang Technological University, Singapore, in 2014, and the B.Sc. degree in communication engineering from Beijing University of Posts and Telecommunications, China, in 2011. He is currently a research scientist with the Beijing Institute for General Artificial Intelligence, Beijing, China.

His research interests include multiagent systems, ad hoc teamwork, deep reinforcement learning, embodied artificial intelligence, and security issues in artificial intelligence. He has published papers in top journals and conferences and has actively served as the reviewer of multiple top conferences for years.
\end{IEEEbiography}
\end{document}